\documentclass[a4paper,11pt]{article}

\usepackage[latin1,utf8]{inputenc}
\usepackage[english]{babel}
\usepackage{amssymb, amsmath,amscd,amsxtra,amsfonts, euscript,latexsym, dsfont, mathrsfs,
 yfonts, longtable,  textcomp, marvosym}
\usepackage{graphics,graphicx}
\usepackage{epsf}
\usepackage{enumerate}
\usepackage{theorem}
\usepackage{fullpage}
\usepackage{dsfont}

\usepackage{
mathptmx,
newcent,
}

\newtheorem{thm}{Theorem}[section]
\newtheorem{defi}[thm]{Definition}
\newtheorem{lemma}[thm]{Lemma}
\newtheorem{prop}[thm]{Proposition}
\newtheorem{ex}[thm]{Example}
\newtheorem{oss}[thm]{Remark}

\newcommand{\boss}{\begin{oss}\rm }
\newcommand{\eoss}{\end{oss}}
\newcommand{\bex}{\begin{ex}\rm }
\newcommand{\eex}{\end{ex}}

\newcommand{\T}{\mathbb{T}}
\newcommand{\rn}{\mathbb{R}^n}
\newcommand{\dimo}{\noindent{\bf Proof\;\, }}
\newcommand{\Min}{\operatorname{Min}}
\newcommand{\res}{\mathop{\hbox{\vrule height 7pt width .5pt depth 0pt
\vrule height .5pt width 6pt depth 0pt}}\nolimits}
\newcommand{\lt}{\{u > t\}}

\newcommand{\qed}{\thinspace\null\nobreak\hfill\hbox{\vbox{\kern-.2pt\hrule
height.2pt depth.2pt\kern-.2pt\kern-.2pt \hbox to2.5mm{\kern-.2pt\vrule
width.4pt \kern-.2pt\raise2.5mm\vbox to.2pt{}\lower0pt\vtop to.2pt{}\hfil
\kern-.2pt \vrule width.4pt\kern-.2pt}\kern-.2pt\kern-.2pt\hrule
height.2pt depth.2pt \kern-.2pt}}\par\medbreak}

\def\ds{\displaystyle}
\def\oF{\overline{F}}
\def\oW{\overline{W}}
\def\A{\mathscr{A}}
\def\R{\mathbb{R}}
\def\Z{\mathbb{Z}}
\def\N{\mathbb{N}}
\def\BV{{\mathrm{BV}}}

\def\wpq#1#2{{\mathrm W}^{#1,#2}}
\def\lp#1{{\mathrm L}^{#1}}
\def\mdiv{\operatorname{div}}
\def\Ind{\operatorname{I}}
\def\Int{\operatorname{Int}}
\def\mod{\operatorname{mod}\;}
\def\cont{{\mathrm{C}}}
\def\Hau{{\mathcal H}}
\def\hone{{\Hau^1}}
\def\one#1{\mathds{1}_{#1}}
\def\vari{{\displaystyle{\mathbf{v}}}}
\def\de{\partial}

\renewcommand{\sectionmark}[1]%
{\markright{\MakeUppercase{\thesection.\#1}}}

\usepackage{hyperref}
\hypersetup{
    backref=true, 
    pagebackref=true,
    plainpages=false,
    bookmarks=true,         
    unicode=false,          
    pdftoolbar=true,        
    pdfmenubar=true,        
    pdffitwindow=true,      
    pdftitle={My title},    
    pdfauthor={Author},     
    pdfsubject={Subject},   
    pdfcreator={Creator},   
    pdfproducer={Producer}, 
    pdfkeywords={keywords}, 
    pdfnewwindow=true,      
    colorlinks=true,       
    linkcolor=blue,          
    citecolor=blue,        
    filecolor=blue,      
    urlcolor=magenta,          
    urlbordercolor=0 1 1
}

\title{A coarea-type formula for the relaxation of a generalized elastica functional}

\author{S. Masnou \footnote{Universit\'e de Lyon, Universit\'e Lyon 1, CNRS (UMR5208), Institut Camille Jordan, 43 boulevard du 11 novembre 1918, F-69622 Villeurbanne-Cedex, France. Email: {\tt masnou@math.univ-lyon1.fr}} \and 
G. Nardi\footnote{Universit\'e Pierre et Marie Curie Paris 6, CNRS UMR 7598, Laboratoire Jacques-Louis Lions, F-75005, Paris,
France. Email: {\tt nardi@ann.jussieu.fr}} }

\date{December 6, 2011}

\begin{document}

\maketitle

\begin{abstract}
We consider the {\em generalized elastica functional} defined on $\lp{1}(\R^2)$ as 
$$F(u)=\left\{\begin{array}{ll}
\ds\int_{\R^2}|\nabla u|(\alpha+\beta|\mdiv
\frac{\nabla u}{|\nabla u|}|^p)\,dx,&\text{if $u\in \cont^2(\R^2),$}\\
+\infty&\text{else},\end{array}\right.$$
where $p>1$, $\alpha>0$, $\beta\geq 0$.
We study the $\lp{1}$-lower semicontinuous envelope $\oF$ of $F$ and we prove that, for any $u\in\BV(\R^2)$, $\oF(u)$ can be represented by a coarea-type formula involving suitable collections of $\wpq{2}{p}$ curves that cover the essential boundaries of the level sets $\{x,\,u(x)> t\}$, $t\in\R$.
\end{abstract}
\section{Introduction}

Being $\BV(\R^2)$ equipped with the strong topology of $\lp{1}$, we consider the functional
$$\begin{array}{rcl}
\mathscr{F}:\BV(\R^2)& \rightarrow& \R\\
u&\mapsto& H(u) + F(u)\end{array}$$
where $H$ is continuous and $\lp{1}$--coercive, i.e.,
$$H(u) \geq \|u\|_1 \quad\quad \forall u\in \BV(\R^2),$$  
and  $F$ is the generalized elastica functional defined on $\lp{1}(\R^2)$ as 
$$F(u)=
\left\{
\begin{array}{ll}
\displaystyle{\int_{\R^2}|\nabla u|(\alpha+\beta|\mdiv
\frac{\nabla u}{|\nabla u|}|^p) dx} & \mbox{if\,} u\in
\cont^2(\R^2),\\
+\infty & \mbox{otherwise.}
\end{array} \right.
$$
with $p> 1$, $\alpha>0$, $\beta\geq 0$, and the convention $|\nabla u||\mdiv \frac{\nabla u}{|\nabla
u|}|^p=0$ whenever $|\nabla u|=0$. When $\beta=0$, the $\lp{1}$-lower semicontinuous envelope $\oF$ is simply (up to the multiplicative constant $\alpha$) the total variation in $\BV$, therefore, to simplify the notations and without loss of generality, we shall assume in the sequel that $\alpha=\beta=1$. The classical Bernoulli-Euler elastica functional associates with any smooth curve $\Gamma\subset\R^2$ its bending energy
$$\int_\Gamma|\kappa_\Gamma|^2d\Hau^1$$
where $\kappa_\Gamma$ is the curvature on $\Gamma$ and $\Hau^1$ the $1$-dimensional Hausdorff measure. The reason why we call $F$ a generalized elastica functional ensues from the fact that, if $u\in\cont^2(\R^2)$, then, by Sard's Lemma, for almost  every $t$, $\partial\{u>t\}\in \cont^2$, $\nabla u/\left|\nabla u\right|(x)$ is orthogonal to $\partial
\{u>t\}$ at every $x$ such that $u(x)=t$, and, by the coarea formula,
$$F(u)=\int_{\R}\left[\int_{\partial \{u>t\}}(1+\left|\mathbf{\kappa}_{\partial\{u>t\}}\right|^p )\,d\mathcal{H}^{1}\right]\mbox{d}t$$

We call {\it $p$-elastica energy} the following map defined on the class of measurable subsets of $\R^2$:
$$E\subset\R^2\mapsto W(E)=\left\{\begin{array}{ll}
\ds\int_{\partial E} [1+|\kappa_{\partial E}|^p] \mbox{d}\Hau^1&\mbox{if }\partial E\in\cont^2,\\
+\infty&\mbox{otherwise}\end{array}\right.
$$
Then,
$$F(u)=\int_{\R}W(\{u>t\})dt\quad\mbox{when }u\in\cont^2(\R^2).$$

A localized variant of $F$ has been introduced in~\cite{MasnouMorel,MM} as a variational model for the inpainting problem in digital image restoration, i.e., the problem of recovering an image known only out of a given domain.  More precisely, it is claimed in~\cite{MasnouMorel,MM} that a reasonable inpainting candidate is a minimizer of this variant of $F$ under suitable boundary constraints. $F$ is also related to a variational model for visual completion arising from a neurogeometric modeling of the visual cortex~\cite{Petitot,CittiSarti}.
As for the numerical approximation of minimizers of $F$, a globally minimizing scheme is
proposed in~\cite{MasnouMorel,MM} for the case $p=1$ while the local minimization
in the case $p>1$ is addressed in~\cite{ChanKangShen} using a fourth-order equation. The case $p>1$ is also tackled in~\cite{CGMP} using a relaxed formulation that involves Euler spirals. A smart  and numerically tractable method to handle the high nonlinearity of the model, actually under a slightly different form, is proposed in~\cite{Ballester}.

The minimization of $F(u)$ in $\BV(\R^2)$ under the constraint that $u$ coincides with a given function out of a given domain has been addressed in~\cite{AM}. The existence of solutions follows from a simple application of the direct method of the calculus of variations. Similarly, proving that the problem 
$$\underset{u\in \BV}{\Min} \mathscr{F}(u)$$
has solutions requires the compactness of minimizing sequences and the lower semicontinuity of $\mathscr{F}$, both with respect to the strong topology of $\lp{1}$. The first property directly follows from the assumptions. Indeed, for every minimizing sequence  $\{u_h\}\subseteq \BV(\R^2)$ with $\underset{h}{\sup}\; \mathscr{F}(u_h)< \infty$, it holds
$$\sup_h\|u_h\|_{\BV}\leq \sup_h\mathscr{F}(u_h)$$
and so, by the compactness theorem for functions of bounded variation, there exists a subsequence converging in $\lp{1}$. The second property, however, does not hold for $F$, as illustrated by the counterexample of Remark \ref{Fsci}, adapted from Bellettini, Dal Maso, and Paolini~\cite{BDP}.

\begin{figure}[h]
\begin{tabular}{c}
\includegraphics[width=15.5cm]{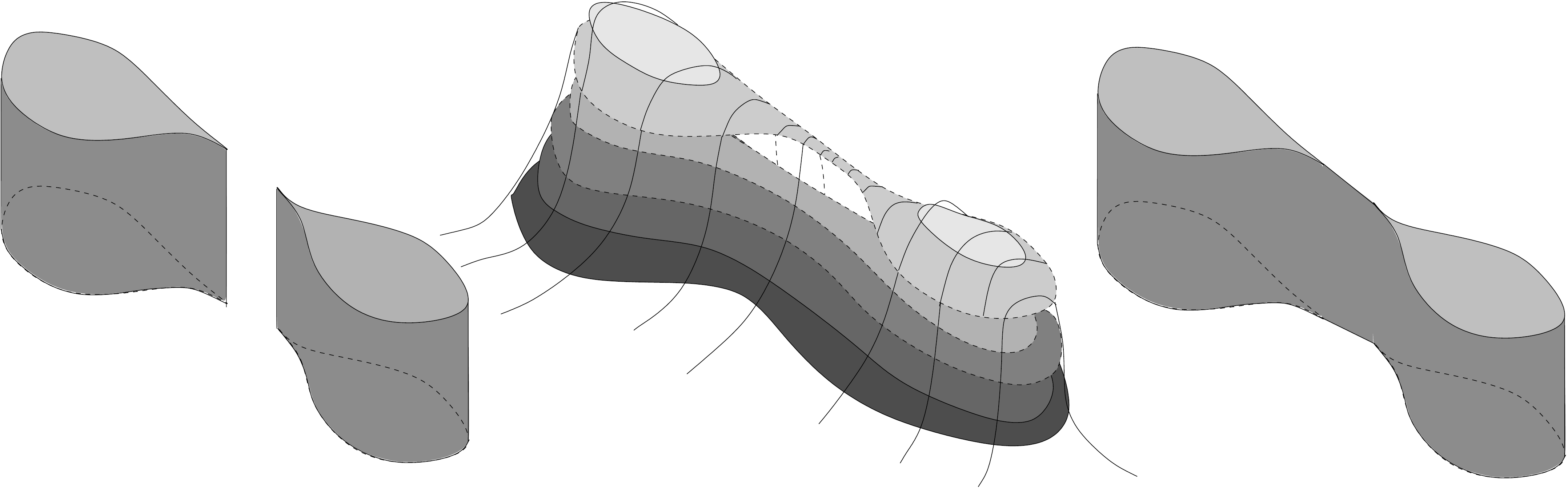}
\end{tabular}
\caption{Left, the characteristic function of a set of finite perimeter. Middle, an approximating smooth function and some of its level sets. Right, the pointwise limit of a sequence of such approximating smooth functions with uniformly bounded elastica energy.}\label{cusp2}
\end{figure}
\begin{oss}[$F$ is not lower semicontinuous in $\lp{1}(\R^2)$]\label{Fsci}

{\rm Let  $E$ be the set of finite perimeter whose characteristic function $u=\mathds{1}_{E}$ is shown in Figure~\ref{cusp2}, left. We can approximate $u$ in $\lp{1}$ by a sequence of smooth functions $\{u_h\} \in
\lp{1}(\R^2)\cap \cont^2_0(\R^2)$, $0\leq u_h\leq 1$, similar to the function partially represented in Figure~\ref{cusp2}, middle. The pointwise limit of such sequence is represented in Figure~\ref{cusp2}, right. By the coarea formula, 
$$F(u_h)=\int_{\R}\left[\int_{\partial \{u_h>t\}}(1+\left|\kappa_{\partial \{u_h>t\}
}\right|^p )d\mathcal{H}^{1}\right]dt$$
and the functions can be designed so that
$$\sup_{h}\, F(u_h)< +\infty.$$
Therefore, 
$$u_h\rightarrow u \mbox{\;in\;}\lp{1}(\R^2)\quad\mbox{and}\quad\liminf_{h\rightarrow +\infty}F(u_h)< +\infty, $$
but, since  $u\notin \cont^2(\R^2)$, we have
$$F(u)=+\infty.$$
Thus the functional $F$ is not lower semicontinuous in $\lp{1}(\R^2)$.}\qed
\end{oss}

As usual for functionals that are not lower semicontinuous~\cite{DalMaso}, we consider the relaxation of $\mathscr{F}$, defined by:
$$\overline{\mathscr{F}}(u)=\inf\left\lbrace
\underset{h\rightarrow\infty}{\liminf}\;\mathscr{F}(u_h) :
u_h\overset{\lp{1}}{\rightarrow}u\right\rbrace .$$
The relaxed functional $\overline{\mathscr{F}}$ is the largest lower semicontinuous functional minoring $\mathscr{F}$ and, in particular,
$$\underset{u\in \BV}{\Min}\; \overline{\mathscr{F}}(u) =  \underset{u\in \BV}{\inf}\; \mathscr{F}(u).$$
In addition, every minimizing sequence of $\mathscr{F}$  has a subsequence converging to a minimum point of $\overline{\mathscr{F}}$ and every minimum point of $\overline{\mathscr{F}}$ is the limit of a minimizing sequence of 
$\mathscr{F}$. More details on the theory of relaxation can be found in~\cite{DalMaso}.
In our case, because of the continuity of $H$, we have
$$\overline{\mathscr{F}} = H + \oF.$$

The existence of minimizers of $\overline{\mathscr{F}}$ by the argument above does not provide much information about $\oF$, about
which only few things are known: it has been proven in~\cite{AM} that
the $N$-dimensional version of $F$ is lower semicontinuous in
$\cont^2(\R^N)$ with respect to the strong topology of $\lp{1}$ when
$N\geq 2$ and $p > N-1$, therefore $\oF(u)=F(u)$ when $u$ is smooth. This constraint on $p$ is weakened in~\cite{LM} where, using results from~\cite{S}, it is shown that $F$ is lower semicontinuous on $\cont^2(\R^2) \cap \lp{1}(\R^2)$ if $N=2$ and $p \geq 1$ or if $N\geq 3$ and $p\geq 2$. Finally, using the same techniques as in~\cite{AM,LM}, and combining with the recent results by Menne~\cite{Menne} on the locality of the mean curvature for integral varifolds, we can conclude that the equality of $F$ and $\oF$ for smooth functions holds in any dimension for any $p\geq 1$. In particular, whenever $N=2$ which is the space dimension for this paper,
$$\oF(u)=F(u)=\int_{\R}\left[\int_{\partial \{u>t\}}(1+\left|\kappa_{\partial \{u>t\}
}\right|^p )d\mathcal{H}^{1}\right]dt=\int_\R W(\{u>t\})\,dt,\quad\mbox{for every }u\in\cont^2(\R^2).$$

We address here a more general question: is there also a coarea-type formula for $\oF(u)$ when $u\in\BV(\R^2)$ has finite relaxed energy ?

\par
This question is obviously related to the relaxation of $W$. The lower semicontinuous envelope of $W$ is defined for every measurable set $E\subset\R^2$ as
$$\oW(E)=\inf\{\liminf_{k\to\infty}\; W(E_k),\; \partial E_k\in\cont^2,\,|E_k\Delta E|\to 0\},$$
where $\Delta$ denotes the symmetric difference operator for sets. We will prove in Proposition~\ref{equiv} that, for any measurable $E\subset\R^2$ such that $\oW(E)<\infty$ and for every $c>0$,
$$\oF(c\one{E})=c\oW(E).$$

The properties of sets with finite relaxed energy $\oW$ have been
extensively studied in~\cite{BDP,BM1,BM}. It is proved in particular
in \cite{BDP,BM1} that, for any $p>1$, $\oW$ can be represented by a
functional depending on systems of curves of class $\wpq{2}{p}$ that
recover and extend the essential boundary of $E$. Another equivalent
representation involving Hutchinson's curvature varifolds is provided
in \cite{BM}. Can these results be used in a straightforward way to
give an explicit expression of $\oF(u)$ for any $u\in\BV(\R^2)$? 

\begin{figure}[h]
\begin{center}
\includegraphics[width=11cm]{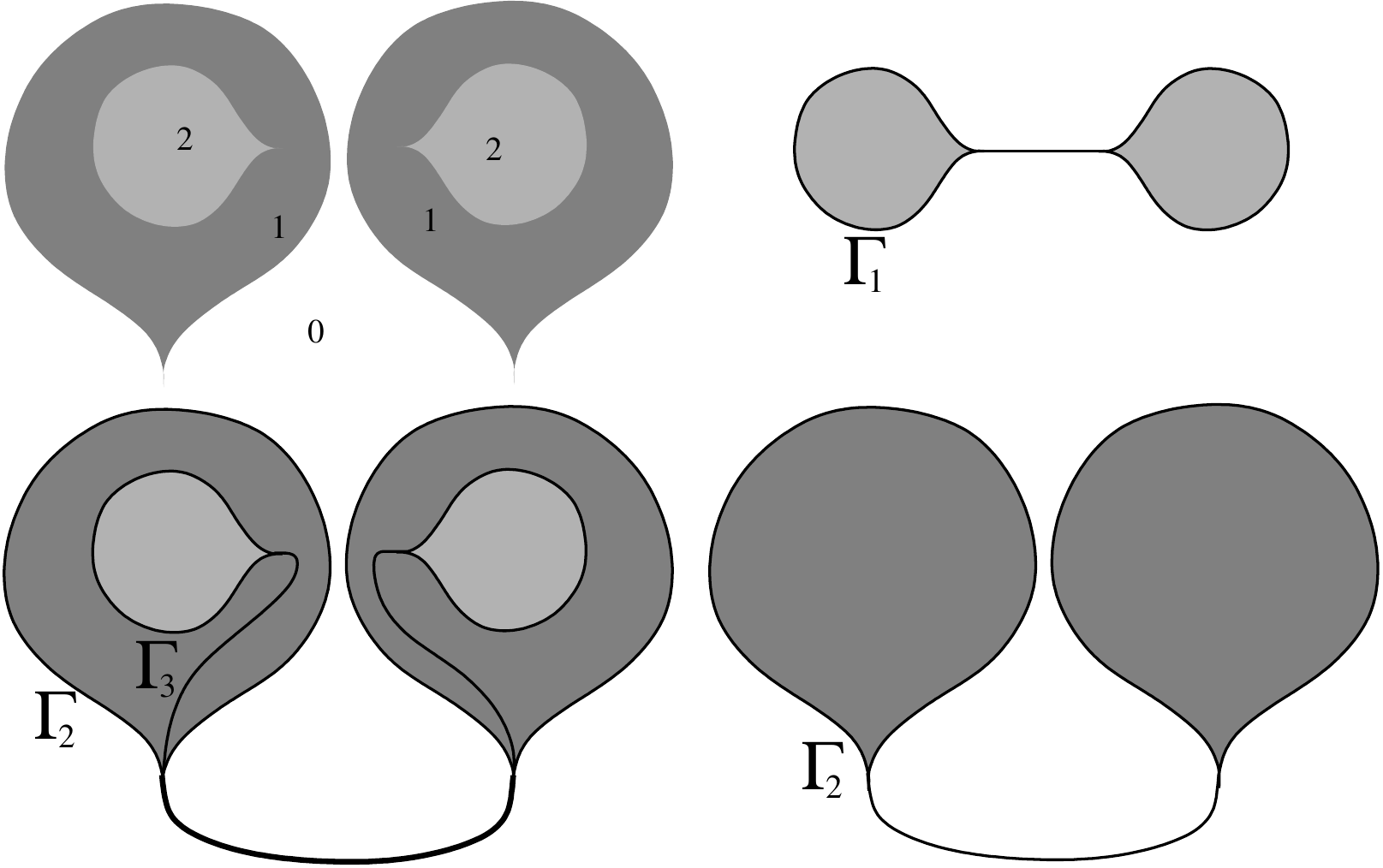}
\caption{The piecewise constant function $u$ on the top left figure has energy $\oF$ stricly greater than $\ds\int_\R\oW(\{u>t\})dt$. Considering an approximating sequence of smooth functions for $\one{\{u\geq 2\}}$, the pointwise singular limit is the bold curve $\Gamma_1$ shown in the top right figure, according to Theorem~8.6 in~\cite{BM1}. The singular limit yielded by the approximation of $\one{\{u\geq 1\}}$ is (approximately) the curve $\Gamma_2$ shown in the bottom right figure. Instead, the smooth approximation of $u$ yields as pointwise singular limit the two curves $\Gamma_2$ and $\Gamma_3$ roughly represented in the bottom left figure. One has $\oF(u)=W(\Gamma_2)+W(\Gamma_3)>W(\Gamma_1)+W(\Gamma_2)=\ds\int_\R\oW(\{u>t\})dt$.}\label{fig:notcoincide}
\end{center}
\end{figure}

A first observation is that, as in the example of Figure~\ref{fig:notcoincide}, $\oF(u)$ and $\ds\int_\R\oW(\partial^*\{u>t\})dt$ do not coincide in general for $u\in\BV(\R^2)$. In the latter integral, $\partial^*\{u>t\}$ denotes the essential boundary of the level set $\{u>t\}$, that has finite perimeter for almost every $t$. Recall that the essential boundary of a set of finite perimeter is the set of points where an approximate tangent exists, see~\cite{AFP}. 

\par
\begin{figure}[h]
\begin{center}
\includegraphics[height=4cm]{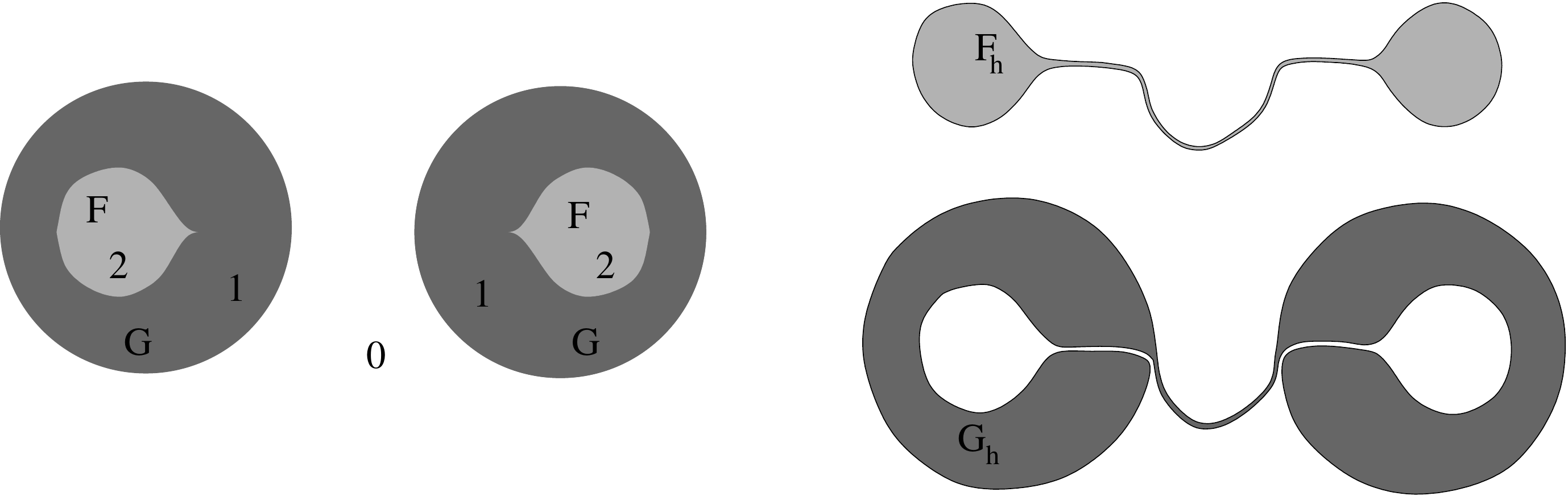}
\caption{The sets $G$ and $F$ can be simultaneously approximated without violating the noncrossing condition}\label{embedded}
\end{center}
\end{figure}

 A more surprising example of a situation where $\oF(u)$ and $\ds\int_\R\oW(\partial^*\{u>t\})dt$ do not coincide is provided in the example below that we shall visit again in Remark~\ref{iii-est-nec}. Let $u=\one{F\cup G}+\one{F}$ with $F, G$ shown in Figure~\ref{embedded}, left. The level set $\{u>0\}$ coincides with $F\cup G$ and is clearly smooth. The boundary of the level  set $\{u>1\}$ has two cusps. The pointwise limit of a sequence of smooth sets that approximate $\{u>1\}$ with convergence of the elastica energy to $\oF(\{u>1\})$ contains the segment joining the two cusps, according to Theorem~8.6 in~\cite{BM1}. By the same theorem, the approximation of $\{u>0\}$ will not contain this segment. Clearly, $\{u>0\}$ and $\{u>1\}$ cannot be contemporaneously approximated using two nested sequences of sets. Yet $\oF(u)$ is finite, as shown using the construction of Figure~\ref{embedded}, right, that we found during discussions with Vicent Caselles and Matteo Novaga. In this construction, the sets $G\setminus F$ and $F$ are approximated using smooth sets that do not intersect. Both pointwise limits contain a "bridge" with multiplicity $2$ that joins both components of $G$.


\par Let us now come back to smooth functions and how their energy simply relates to the energy of their level sets. Given $u\in\BV(\R^2)$ such that $\oF(u)<\infty$, one considers a
sequence of smooth functions $(u_h)$ converging to $u$ in
$\lp{1}(\R^2)$ and such that $F(u_h)\to \oF(u)$ as $h\to\infty$. Possibly extracting a subsequence, one can assume that for almost every $t$, $\{u_h>t\}$ converges to $\{u>t\}$ in measure, i.e. $|\{u_h>t\}\Delta\{u>t\}|\to 0$. In addition, by Fatou's
Lemma,
$$\int_\R\liminf_{h\to\infty}\,W(\{u_h>t\})dt\leq\liminf_{h\to\infty}\int_\R W(\{u_h>t\})dt=\oF(u)$$
therefore $\liminf_{h\to\infty}W(\{u_h>t\})<\infty$ is finite for almost every $t$. It follows that, for almost every $t$, the sequence of $1$-dimensional varifolds with unit multiplicity $\vari(\partial\{u_h>t\},1)$ has uniformly bounded mass, and uniformly bounded curvature in $\lp{p}$. By the properties of varifolds~\cite{Si} and the stability of absolute continuity (see Example 2.36 in~\cite{AFP}), there exists a subsequence $\vari(\partial\{u_{h_k}>t\},1)$ \textit{depending on $t$} and a limit integral $1$-varifold $V_t$ such that 
$$\int_{\R^2}(1+|\kappa_{V_t}|^p)d\|V_t\|\leq \liminf_{h\to\infty}\,W(\{u_h>t\})$$
In addition, one can prove~\cite{AM} that the support $M_t$ of $V_t$ contains $\partial^*\{u>t\}$ for almost every $t$. Furthermore, if $u$ is \textit{smooth}, $\kappa_{V_t}$ coincides almost everywhere on $\de\{u>t\}$ with $\kappa_{\de\{u>t\}}$ thus
$$W(\{u>t\})=\int_{\de\{u>t\}}(1+|\kappa_{\de\{u>t\}}|^p)d\Hau^1\leq\int_{\R^2}(1+|\kappa_{V_t}|^p)d\|V_t\|\leq\liminf_{h\to\infty}\,W(\{u_h>t\}).$$
By a simple integration and the coarea formula, it follows that for any smooth function $u$, $F(u)=\oF(u)$. The argument above is exactly the $2$-dimensional version of the more general proof proposed in~\cite{AM,LM} to prove the equality of $F$ and $\oF$ on smooth functions in any dimension for various ranges of values of $p$.

Can a similar argument be used if $u$ is \textit{unsmooth}? A tentative strategy could be the following:
\begin{enumerate}
\item show, if possible, that the limit varifolds $V_t$ built above are nested, i.e. $\operatorname{int}V_t\subset\operatorname{int}V_{t'}$ if $t>t'$, where $\operatorname{int}V_t$ denotes the set enclosed (in the measure-theoretic sense) by the support of $V_t$. Again, observe that $\ds\int_{\R^2}(1+|\kappa_{V_t}|^p)d\|V_t\|\leq \liminf_{h\to\infty}\,W(\{u_h>t\})$.
\item using the results of~\cite{BDP}, build a sequence of sets $E_h^t$ (for a suitable dense set of values $t$) such that $\partial E_h^t\to M_t$ (being $M_t$ the support of $V_t$) and $W(E_h^t)\to\ds\int_{\R^2}(1+|\kappa_{V_t}|^p)d\|V_t\|$. The varifolds $V_t$ being nested, one could actually build $E_h^t$ so that $E_h^t\subset E_h^{t'}$ if $t>t'$.
\item by a suitable smoothing of the sets $E_h^t$, build a smooth function $\tilde u_h$ such that $F(\tilde u_h)\leq \ds\int_\R W(E_h^t)dt+\frac 1 h$. 
\item passing to the limit, possibly using a subsequence, show that $\tilde u_h$ tends to $u$ in $\lp{1}$ and using the lower semicontinuity of $\oF$, conclude that
$$\oF(u)=\int_\R\int_{\R^2}(1+|\kappa_{V_t}|^p)d\|V_t\|\,dt$$
\end{enumerate}

This strategy has however a major difficulty: the fact that the limit varifolds are nested is not clear at all. It would be an easy consequence of the existence of a subsequence $(u_{h_k})$ such that the varifolds $\vari(\partial\{u_{h_k}>t\},1)$ converge to $V_t$ for almost every $t$. But such subsequence {\it may not exist} in general as shown by the counterexample below due to G. Savar{\'e}~\cite{savare}.

\begin{ex}[Savar\'e~\cite{savare}]\label{savare}{\rm  Let us design a sequence of functions $\{\tilde{u}_n\}\subset \cont^{0}([0,1]^2)$ with smooth level lines $\{\tilde u_n=t\}$ satisfying 
$$\sup_n \int_{\R}  \int_{\partial \{\tilde{u}_n(x)>t\}\cap (0,1)^2}(1+\left|\kappa_{\partial \{\tilde{u}_n(x)>t\}\cap (0,1)^2}\right|^p ) \,d\mathcal{H}^{1} \,dt <  \infty,$$
but such that there exists no subsequence $(t\mapsto\vari(\partial\{\tilde u_{n_k}>t\},1))_k$ converging for almost every $t$ to a varifold $V_t$. 
\par\noindent  Consider  the following 2-periodic  function on $\R$ :
$$u(x)=
\left\{
\begin{array}{rl}
1 & \mbox{\;if \,} x\in ]0,1/2[\cup]1,2[, \\
-1 & \mbox{\;if \,}x\in ]1/2,1[
\end{array} \right.
$$
\noindent Define $u_n(x)= u(2^n x)$ and consider 
$$U(x)=\int_0^x u(s) \mbox{d}s\, , \quad \quad U_n(x)=\int_0^x u_n(s) \mbox{d}s.$$

\begin{figure}[h]
\begin{center}
\includegraphics[totalheight=6cm]{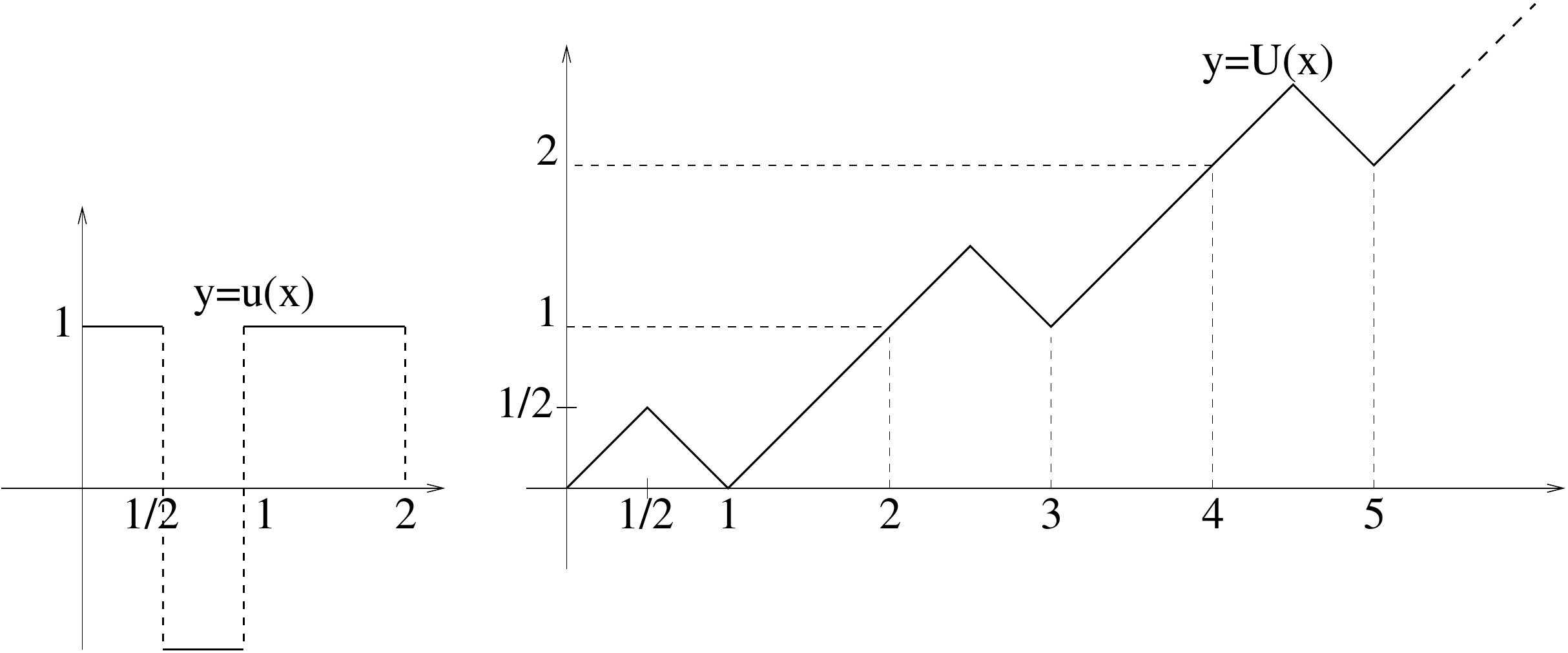}
\caption{Graphs of the functions $u$ and $U$}
\end{center}
\end{figure}

\noindent The sequence $\{U_n\}$ is equilipschitz, because $|U_n'(x)|=1$, and  
$$U_n(x)=\int_0^x u_n(s) \mbox{d}s = \dfrac{1}{2^n}\int_0^{2^n x} u(s) \mbox{d}s = \dfrac{1}{2^n} U(2^n x).$$
Thus
$$\{U_n(x)=t\}=\{U(2^n x)= 2^n t\}\quad \forall t\in \R.$$
Then
$$\mathcal{H}^0\left( \{U(x)=\sigma \}\right) =
\left\{
\begin{array}{rl}
3 & \mbox{\,\;if \,} \sigma\in ]k,k+1/2[, \\
1 & \mbox{\,\;if \,} \sigma\in ]k+1/2,k+1[
\end{array} \right.\quad\forall k\in\N,
$$
and so, for almost every $t$,
$$\mathcal{H}^0\left( \{U_n(x)=t \}\right) =
\left\{
\begin{array}{rl}
3 & \mbox{\,\;if \,} t\in ] \frac{2k}{2^{n+1}}, \frac{2k+1}{2^{n+1}} [,  \\
1 & \mbox{\,\;if \,} t\in ] \frac{2k+1}{2^{n+1}}, \frac{2k+2}{2^{n+1}} [
\end{array} \right.\quad\forall k\in\N.
$$
\par\noindent We now define the sequence of functions 
$$\tilde{u}_n(x_1, x_2)=U_n(x_1) \quad \quad x=(x_1,x_2)\in [0,1]^2$$
and consider the sequence of varifolds associated with its level lines
$$V_t^n = {\bf v}(\partial\{\tilde{u}_n(x)>t\}\cap (0,1)^2, 1).$$
Remark that $U_n(1)=1/2$ for every $n$ and, by the coarea formula, it is easy to check that
$$\sup_n \int_{\R} \int_{\partial \{\tilde{u}_n(x)>t\}\cap (0,1)^2}(1+\left|\kappa_{\partial \{\tilde{u}_n(x)>t\}\cap (0,1)^2}\right|^p )  \,d\mathcal{H}^{1} \,dt< 3/2.$$
Nevertheless we can show that there exists no subsequence $\{V^n_t\}$ (not relabelled) and no family $(V_t)_{t\in \R}$ of varifolds such that $V^n_t \rightharpoonup V_t$ in the sense of varifolds for almost every $t$.

\par For that we can consider the sequence $\{\mu_{V^n_t}\}$ of the weight measures of the varifolds $V^n_t$ defined as
$$\mu_{V^n_t} = \mathcal{H}^1(\partial \{\tilde{u}_n(x)>t\}\cap (0,1)^2)$$
therefore
$$\mu_{V^n_t} = \mathcal{H}^0\res \{U_n(x)=t\}.$$
By contradiction we suppose that there exists a subsequence $(\mu_{V^n_t})$ (not relabelled) and a family of limit measures $(\mu(t))_{t\in \R}$ such that 
$$\langle \mu_{V^n_t},\varphi \rangle=\langle\mathcal{H}^0\res \{U_n(x)=t\}, \varphi \rangle \rightarrow \langle\mu(t), \varphi\rangle \quad \text{for a.e.}\; t \quad \forall \varphi \in \cont^1_0([0,1]).$$

\par Define
$$f_n(t) = \mathcal{H}^0\res\{U_n(x)=t\}([0,1]).$$
Since $\{f_n\}$ is uniformly bounded  we get  
$$f_n(t) \rightarrow \mu(t)([0,1]) \quad\text{for a.e. $t$}$$
 and so $\mu(t)([0,1])$ is bounded. Then, by the Dominated Convergence Theorem,  we get 
$$f_n(t) \rightarrow \mu(t)([0,1]) \quad \text{in} \quad \lp{2}(]0,2[).$$ 
\par On the other hand $f_n(t)=f(2^{n+1}t)$ where $f$ is the 2-periodic function on $\R$ defined as
$$f(t)=
\left\{
\begin{array}{rl}
3 & \mbox{\;if \,} t\in ]0,1[ \\
1 & \mbox{\;if \,} t\in ]1,2[
\end{array}  \right.
$$
and so, by the Riemann-Lebesgue Theorem, it ensues that $f_n \rightharpoonup 2$ weakly in $\lp{2}$, therefore, by the strong convergence in $\lp{2}$ established above, we should have $\mu(t)([0,1])=2$. But $2$ cannot be the strong limit of $\{f_n\}$ thus there exists no subsequence of $\{\mathcal{H}^0\res \{U_n(x)=t\}\}$ converging for almost every $t$ to a limit measure $\mu(t)$.\qed}
\end{ex}

\par It is now clear that, given a sequence of smooth functions $(u_h)$ converging in $\lp{1}$ to $u\in\BV$ with uniformly bounded generalized elastica energy for some $p>1$, we cannot expect in general the convergence of a subsequence $(\de\{u_{h_k}>t\})$ to a limit system of curves $\Gamma_t$ for almost every $t$. Instead {\it we will prove in this paper that we can find a countable and dense set of values $I$ such that $\de\{u_{h_k}>t\}\to \Gamma_t$ for every $t\in I$}. Then, for almost every remaining value $t\in\R\setminus I$, a limit system of curves can be built as the limit of a sequence $(V_{t_n})$, $t_n\in I$, $t_n\to t$. This "dense" diagonal extraction is obtained by generalizing the approach used in~\cite{MM} to study the same functional $\oF$ on a restricted class $\mathcal{S}$ of functions $u\in \BV(\R^2)$ such that, for a given $\Omega$ with smooth boundary $\partial\Omega\in\cont^\infty$,
\begin{itemize}
 \item $u=U_0$ on $\R^2\setminus \Omega$ with $U_0$ analytic on $\R^2$, $F(U_0)<+\infty$ and the level lines of $U_0$ in $\Omega$ satisfy a few regularity assumptions;
 \item for $\mathcal{L}^1$-a.e. $t$ the restriction to $\Omega$ of $\partial^* \{u>t\}$ coincides, up to a $\mathcal{H}^1$-negligible set, with the trace of finitely many curves of class $\wpq{2}{p}$ with each of them joining two points on $\partial\Omega$ and smoothly connected to a level line of $U_0$ out of $\Omega$ (see~\cite{MM} for details).
\end{itemize}
In this paper, the boundary constraints are dropped, in particular the level lines are no more constrained to link two points on a given boundary. Dropping this constraint raises a few technical difficulties, in particular in some situations of accumulation that will be detailed later on. The "dense diagonal" convergence technique mentioned above was used in~\cite{MM} to build a set of curves that cover the level lines of a minimizer of $\oF$ in the class ${\cal S}$. We shall here extend this convergence technique to obtain, for every $u\in\BV(\R^2)$ with $\oF(u)<+\infty$, a coarea-type representation formula for $\oF(u)$ using the $p$-elastica energies of curves that cover the essential boundaries of the level sets of $u$. To be more precise, we will associate with each $u\in\BV(\R^2)$ having finite energy a class of functions $\A(u)$ defined as follows: $\Phi\in\A(u)$ whenever $\Phi:t\in\R\mapsto \{\gamma_t^1,\cdots,\gamma_t^N\}$ with $N$ depending on $t$ and $\{\gamma_t^1,\
 cdots,\gamma_t^N\}$ is a finite collection of curves of class $\wpq{2}{p}$ without crossing (but possibly with tangential self-contacts) that satisfy nesting compatibility constraints. Defining the functional 
$$\Phi\in \mathscr{A}(u) \mapsto G(\Phi)=\int_{\R} W(\Phi(t)) \mbox{d}t,$$
we will show in Theorem~\ref{principale} that for every  $u\in \BV(\R^2)$ with  $\oF(u) < \infty$ there holds
$$\oF(u) = \underset{\Phi \in \mathscr{A}(u)} {\Min} G(\Phi),$$
which is the desired coarea-type representation formula.
The main difficulty in the proof is to handle properly the situations of accumulation that possibly occur for the graphs of approximating functions.

Let us conclude this introduction with a short comment about the problem in higher dimensions. Is there a similar decomposition of $\oF$ using suitable covering of level hypersurfaces? It is an open problem to our knowledge and, anyway, the solution needs not involve finite collections of hypersurfaces at each level. Indeed, an example due to Brakke~\cite{Br} consists of an integral $2$-varifold in $\R^3$ with uniformly bounded mean curvature and such that, at no point of a set with positive measure, the varifold's support can be represented as the graph of a multi-function. Even the control of the whole second fundamental form is not enough: it is shown in \cite[Thms 3.3 and 3.4, Example 5.9]{AGP}) that if $p>N\geq 3$, the limit varifold of a sequence of  smooth  boundaries with equibounded $\lp{p}$-norm of the second fundamental form needs not be representable as a finite union of manifolds of class $\wpq{2}{p}$. In a companion paper~\cite{MasnouNardi2}, we propose instead 
 a completely different strategy based on varifolds associated with gradient Young measures.

The plan of the paper is as follows: in Section~\ref{notation} we introduce a few notations, and we define and discuss the class $\mathscr{A}(u)$ mentioned above. We prove in Section~\ref{G} that the minimum problem for $G$ has at least a solution in $\mathscr{A}(u)$, and in Section~\ref{formula} we show the characterization formula for $\oF$. We further illustrate in Section~\ref{oFoW} the connection between $\oF$ and $\oW$. Lastly, in Section~\ref{omega}, we analyze the generalized elastica functional localized on a domain $\Omega\subset \R^2$.

\section{Notations and preliminaries}\label{notation}

Throughout the paper, $p>1$ is a real number, $\mathcal{L}^n$ the Lebesgue measure on $\R^n$,  ${\mathcal H}^k$ the $k$-dimensional Hausdorff measure, and $\cont^r$, $\lp{p}$, $\wpq{m}{p}$, $\BV$  the usual function spaces. For any $E \subseteq \rn$ we will also denote $|E|=\mathcal{L}^n(E)$. The topological boundary of $E$ is denoted as $\partial E$ and, if $E$ has finite perimeter~\cite{AFP}, $\partial^* E$ is its essential boundary, i.e. $\partial^*E = \R^n\setminus (E_0\cup E_1)$ where, for every $t\in [0,1]$,
$$E_t= \left\lbrace x\in \R^n : \underset{r \rightarrow 0}{\lim} \dfrac{|E\cap B(x,r)|}{|B(x,r)|} = t\right\rbrace .$$
In addition, by Federer Theorem~\cite[Thm 3.61]{AFP}, $\partial^*E$ coincides, up to a ${\mathcal H}^{N-1}$-negligible set, with the set of points where the inner normal $\ds\frac{D\one{E}}{|D\one{E}|}$ exists.

If the topological boundary $\partial E$ can be viewed, locally, as the graph of a function of class $\cont^r$ (resp. $W^{m,p}$), we write $\partial E \in \cont^ r$ (resp. $W^{m,p}$).

\par Unless specified, we now focus on two-dimensional sets. We first recall the definition of the index of a point with respect to a plane curve $\gamma$~\cite{rudin-real-complex}.

\begin{defi}
 Let $\gamma : [0,1] \rightarrow \mathbb{C}$ be a $\cont^1$-curve with support $(\gamma)=\gamma([0,1])$. The index of a point $p\in \mathbb{C}\setminus(\gamma)$ with respect to $\gamma$ is defined by
$$\Ind(p, \gamma) = \dfrac{1}{2\pi i} \int_{\gamma}\dfrac{\,dz}{z-p}\cdot$$

\end{defi}
In the sequel, we denote as $k_{\gamma}(s)=\gamma''(s)$ the curvature vector at a point $\gamma(s)$ of a curve $\wpq{2}{p}$-curve $\gamma$ parameterized with arc-length $s$, i.e. at constant unit velocity. We shall use the following convenient and classical lemma, whose proof is recalled.

\begin{lemma}\label{mon}
 Let $\mathcal{F} = \left( X_t \right)_{t\in \R}$  be a monotone family of sets, $X_t \subseteq \rn$ for all $t$. Then, there exists an at most countable set $D \subseteq \R$ such that for every compact set $K\subset \rn$
$$\underset{s\rightarrow t}{\lim}\;  |(X_s \Delta X_t)\cap K| = 0 \quad  \forall t\in \R\setminus D.$$
We call $D$ the set of discontinuities of $\mathcal{F}.$
\end{lemma}

\dimo The family $\mathcal{F}$ is monotone so the function 
$$t \mapsto | X_t \cap K|$$
is monotone for every compact set $K$ and it has at most  countably many  discontinuity points whose collection is denoted as $D$. Then for every $t\in \R\setminus D$ we have  $| X_s \cap K| \rightarrow | X_t \cap K| $ as $s \rightarrow t$ and because of the monotonicity of $\mathcal{F}$ we get 
$$\underset{s\rightarrow t}{\lim}\;  |(X_s \Delta X_t)\cap K| = 0 \quad  \forall t\in \R\setminus D.$$
\qed

\noindent Following~\cite{BDP}, we now define the notion of system of curves of class $\wpq{2}{p}$.

\begin{defi}\label{sistemi}
 By a system of curves of class $\wpq{2}{p}$ we mean a finite family $\Gamma=\{\gamma_1,\cdots,\gamma_N\}$ of closed curves of class $\wpq{2}{p}$ (thus $\cont^1$) admitting a parameterization  (still denoted by  $\gamma_i$) $\gamma_i\in \wpq{2}{p}\left([0,1],\R^2 \right)$  with constant velocity.
Moreover, every curve of  $\Gamma$ can have  tangential self-contacts but without crossing and two curves of $\Gamma$ can have  tangential contacts but without crossing. In particular, $\gamma_i'(t_1)$ and $\gamma_j'(t_2)$ are parallel whenever $\gamma_i(t_1)=\gamma_j(t_2)$ for some $i,j\in\{1,...,N\}$ and $t_1,t_2\in[0,1]$.
\par The trace  $(\Gamma)$ of $\Gamma$ is the union of the traces $(\gamma_i)$. We define the interior of the system $\Gamma$ as
$$ \Int(\Gamma) = \{x \in \R^2 \setminus (\Gamma) : \Ind(x, \Gamma) = 1 \; \mod 2\}$$
where $\Ind(x, \Gamma) = \sum_{i=1}^N \Ind(x, \gamma_i)$.\\
The multiplicity  function $\theta_\Gamma$ of $\Gamma$ is
$$\theta_{\Gamma} :  (\Gamma) \rightarrow \mathbb{N} \quad \quad \theta(z)=\sharp \{\Gamma^{-1}(z)\},$$
where $\sharp$ is the counting measure. \\
If the system of curves is the boundary of a set $E$ with $\partial E\in \cont^2$, we simply denote it as $\partial E$.
\end{defi}

\begin{oss}{\rm Remark that, by previous definition, every $|\gamma_i'(t)|$ is constant for every  $t\in[0,1]$ so the arc-length parameter 
is given by $s(t)=t L_i$ where $L_i$ in the length of $\gamma_i$. Denoting by $\tilde{\gamma}_i$ 
the curve parameterized with respect to the arc-length parameter we have  
$$s\in[0, L_i] \,,\,\, \tilde{\gamma}_i(s)= \gamma_i(s/L_i)\,,\,\, \tilde{\gamma}_i''(s)=\dfrac{\gamma_i''(s)}{L_i^2}.$$
Now, the curvature {\bf k} as a functions of $s$, verifies
$${\bf k} = \tilde{\gamma}_i''(s)$$
which implies
$$ \int_{0}^{L_i} \left( 1+|\tilde{\gamma}_{i}''(s)|^p\right)\mbox{d}s=\int_{0}^{L_i} \left( 1+|k|^p\right)\mbox{d}s= \int_{0}^{1} \left( |\gamma_{i}'(t)|+L_i^{1-2p}|\gamma_{i}''(t)|^p\right)\mbox{d}t.$$
Then, the condition $\gamma_i \in \wpq{2}{p}\left([0,1],\R^2 \right)$ implies that $\tilde{\gamma}_i\in \wpq{2}{p}\left([0,L_i],\R^2 \right)$ 
and, for simplicity, in the sequel we denote  by $\gamma_i$  the curve parameterized with respect to the arc-length parameter.}
\end{oss}

The $p$-elastica energy of a system $\Gamma$ of curves of class $\wpq{2}{p}$ is defined as
$$W(\Gamma)=\sum_{i=1}^N W(\gamma_i)=\sum_{i=1}^N \int_{(\gamma_i)} \left( 1+|{\bf k}_{\gamma_{i}}|^p\right)\mbox{d}\mathcal{H}^1.$$

We will use several times the following result that combines Lemma 3.1 in~\cite{BDP} and Proposition 6.1 in~\cite{BM1}: if a bounded open set $E\subset\R^2$ is such that $\oW(E)<\infty$ then
$$\overline{W}(E)=\min\;\{W(\Gamma): \Gamma \in \mathcal{A}(E)\}$$
where $\mathcal{A}(E)$ denotes the class of all {\it finite} systems of curves $\Gamma$ such that $(\Gamma)\supseteq \partial E$ and $|E \Delta \Int(\Gamma)| = 0$.
\begin{defi}[Convergence of systems of curves]
 Let $(\Gamma_h)_{h\in\N}=(\{\gamma_1^h,...,\gamma_{N(n)}^h\})_h$ be a sequence of system of curves of class $\wpq{2}{p}$. We say that $(\Gamma_h)$ converges weakly in $\wpq{2}{p}$ to $\Gamma=\{\gamma_1,...,\gamma_M\}$ if
\begin{itemize}
 \item[(i)] $N(h)=M$ for $h$ large enough;
 \item[(ii)] $\gamma_i^h$ converges weakly in $\wpq{2}{p}$ to $\gamma_i$ for every $i\in\{1,\cdots,M\}$.
\end{itemize}
\end{defi}

\begin{defi}\label{sistlim}
We say that $\Gamma$ is a limit system of curves of class $\wpq{2}{p}$ if $\Gamma$ is the weak limit of a sequence $(\Gamma_h)$ of boundaries of bounded open sets with $\wpq{2}{p}$ parameterizations.
\end{defi}


\begin{defi}\label{A}
Let $\mathscr{A}$ denote the class of functions
$$\Phi : t\in\R \rightarrow \Phi(t) $$
where for almost every $t\in \R$, $\Phi(t)=\{\gamma_t^1,...,\gamma_t^N\}$ is a limit system of curves of class $\wpq{2}{p}$ and such that,  for almost every $\underline{t}, \overline{t}\in \R$, $\underline{t}< \overline{t}$, the following conditions are satisfied: 
\begin{itemize}
 \item[(i)]  $\Phi(\underline{t})$ and $\Phi(\overline{t})$ do not cross but may intersect tangentially;
 \item[(ii)] $\Int(\Phi(\overline{t})) \subseteq \Int(\Phi(\underline{t}))$ (pointwisely);
 \item[(iii)] if, for some $i$, $\mathcal{H}^1\left( (\gamma_{\overline{t}}^i) \setminus \overline{\Int(\Phi(\underline{t}))}\right) \neq 0$ then 
$$\mathcal{H}^1\left( [(\gamma_{\overline{t}}^i) \setminus \overline{\Int(\Phi(\underline{t}))}] \setminus (\Phi(\underline{t}))\right) = 0.$$ 
\end{itemize} 
\end{defi}

\boss\label{remarkiii}
 One may remark that, from condition $(ii)$ of Definition~\ref{A}, for every curve  $\gamma\in\Phi(\underline{t})$
$$ (\gamma)  \cap  \Int(\Phi(\overline{t}))  =\emptyset.$$
In fact if $x\in (\gamma)  \cap  \Int(\Phi(\overline{t}))$ then $x\in \Int(\Phi(\overline{t}))$ and $x\notin  \Int(\Phi(\underline{t}))$ which gives a contradiction with condition $(ii)$.
\eoss

Following~\cite{MM},  we introduce a convenient notion of convergence in $\mathscr{A}$: 

\begin{defi}[Convergence in $\mathscr{A}$]\label{convA}
We say that $\Phi_h$ converges to $\Phi$ in $\mathscr{A}$, and we denote $\Phi_h \overset{\mathscr{A}}{\longrightarrow}\Phi$,  if 
\begin{itemize}
 \item[(i)] for each dyadic interval $[k_N2^{-N}, (k_N+1)2^{-N})$, $N\geq 1$, $k_N\in\Z$, there exists  a point $t_{N, k_N}$ in the interval such that $\Phi_h(t_{N, k_N})$
converges to $\Phi(t_{N, k_N})$ weakly in $\wpq{2}{p}$ as $h\rightarrow \infty$;
\item[(ii)] for almost every  $t\in \R$, there exists a sequence $\{t_{N,k_N}\}$ such that   $t_{N,k_N}\rightarrow t$ and  $\Phi(t)$ is the
weak $\wpq{2}{p}$ limit  of   $\{\Phi(t_{N, k_N})\}$ as $N\rightarrow \infty$.
\end{itemize}
\end{defi} 
It follows from this definition that, if $\Phi_h \overset{\mathscr{A}}{\longrightarrow}\Phi$, there exists for almost every  $t\in \R$ a
sequence $(h_N, k_N)$  such that 
$$ t_{N,k_N}\rightarrow t,\quad h_N\rightarrow \infty$$
and  
$$\Phi_{h_N}(t_{N,k_N})\overset{\wpq{2}{p}}{\rightharpoonup}  \Phi(t)$$
therefore
$$\Phi_{h_N}(t_{N,k_N})\overset{\cont^1}{\longrightarrow}  \Phi(t).$$
 
In the following definition, we associate with any function $u$ of bounded variation in the plane the class of all functions in $\mathscr{A}$ that realize a nested covering of the essential boundaries of the level sets of  $u$, i.e. a covering of its level lines.

\begin{defi}[The class $\mathscr{A}(u)$]\label{defiA}
Let $u\in \BV(\R^2)$. We define $\mathscr{A}(u)$ as the set of functions
$\Phi \in \mathscr{A} $ such that, for almost every $t\in \R$, we have 
$$(\Phi(t))\supseteq \partial^* \lt \quad \mbox{ (up to a $\mathcal{H}^1$-negligible  set)}$$
 and 
$$\lt = \Int(\Phi(t)) \quad \mbox{ (up to a ${\cal L}^2$- negligible  set)}.$$
In particular, if $u\in \cont^2(\R^2)$, we will denote as $\Phi[u]$ the function of  $\mathscr{A}(u)$ defined as
$$t \mapsto \partial \{u>t\}.$$
\end{defi}
{\boss We will prove in Theorem~\ref{minG} that, whenever $u\in\BV(\R^2)$ is such that $\oF(u)<\infty$, then $\mathscr{A}(u)\not=\emptyset$.
\eoss}

Conditions in Definitions~\ref{A} and~\ref{defiA} ensure that any $\Phi\in\A(u)$ is a nested covering of the level lines of $u$. In particular, condition $(iii)$ of 
Definition~\ref{A} ensures that the nesting property is also satisfied wherever concentration occurs, typically on the ghost concentration segment of Figure~\ref{fig:notcoincide}, right. The following examples show the necessity of condition~$(iii)$.
\begin{oss}[Condition $\ref{A}(iii)$ is necessary]\label{iii-est-nec}
{\rm 
Let $u=\mathds{1}_{E\cup F} + \mathds{1}_F$ and $v=\mathds{1}_{G\cup F} + \mathds{1}_F$ with $E,F,G$ like in Figure~\ref{condintro}-A.
\begin{figure}[h]
\begin{center}
\includegraphics[height=8cm]{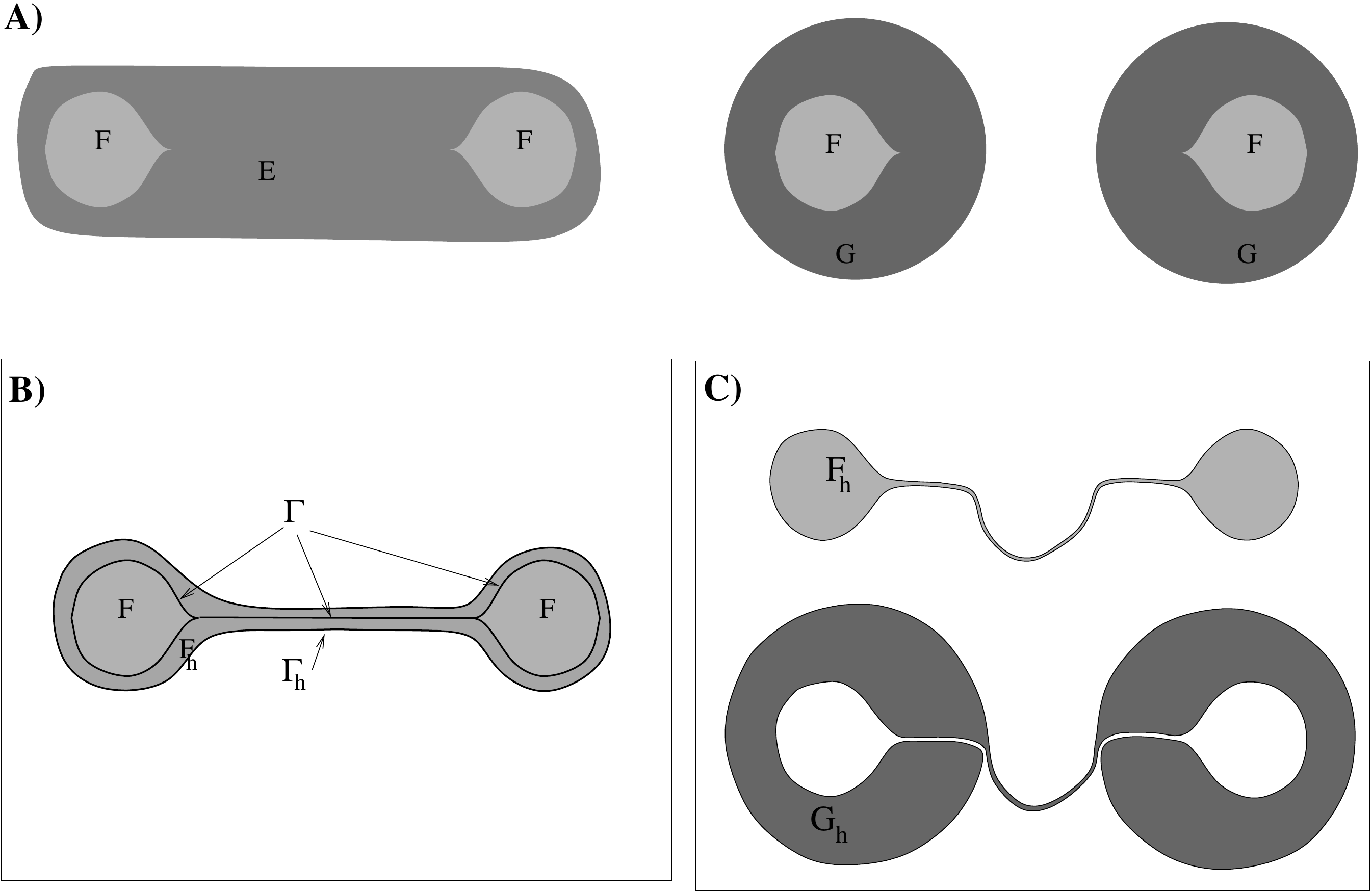}
\caption{Situations of accumulation}\label{condintro}
\end{center}
\end{figure}
Using the same kind of approximation as in Figure~\ref{cusp2} one easily sees that $u$ has finite energy $\oF$. Sequences of smooth functions $(u_h)$ that converge to
 $u$ in $\lp{1}$ and such that $F(u_h) \rightarrow \oF(u)$ have level sets similar to $F_h$ in Figure~\ref{condintro}-B for every level between 1 and 2. 
\par The boundary $\Gamma_h$ of $F_h$ converges to the limit curve $\Gamma$ that contains $\partial F$ plus a ghost segment that corresponds to the concentration in the limit of the middle tube. Clearly,  for $h$ large enough, the middle tube is contained in sets that approximate $E$.

The situation is different for $v$ which also has finite relaxed energy, as was explained in the introduction. To build the approximating sequence with uniformly bounded energy, it is necessary to approximate $G\setminus F$ in such a way that a ``corridor'' be created between both disks so that both components of $F$ can be approximated with a sequence of connected sets having bounded energy. This is illustrated in Figure~\ref{condintro}-C: the bottom set is smooth and connected. As the width of all thin gray and white zones goes to $0$, the set converges in measure to $G\setminus F$. This clearly justifies the need for condition $\ref{A}(iii)$ since the simultaneous approximation of $F$ and $G$ without boundary crossing requires building also for $G\setminus F$ a ghost part that encloses the ghost part arising from the approximation of $F$.}
\end{oss}

\begin{oss}[Exemplification of Definition~\ref{A}]{\rm Let us analyze on some examples the geometric meaning of Definition \ref{A}.
\begin{itemize}
\item{\it Example I}\;  Let $u = \mathds{1}_E + \mathds{1}_F$ with $E,F$ like in Figure~
\ref{wall1}.

\begin{figure}[!h]
\begin{center}
\includegraphics[height=3cm]{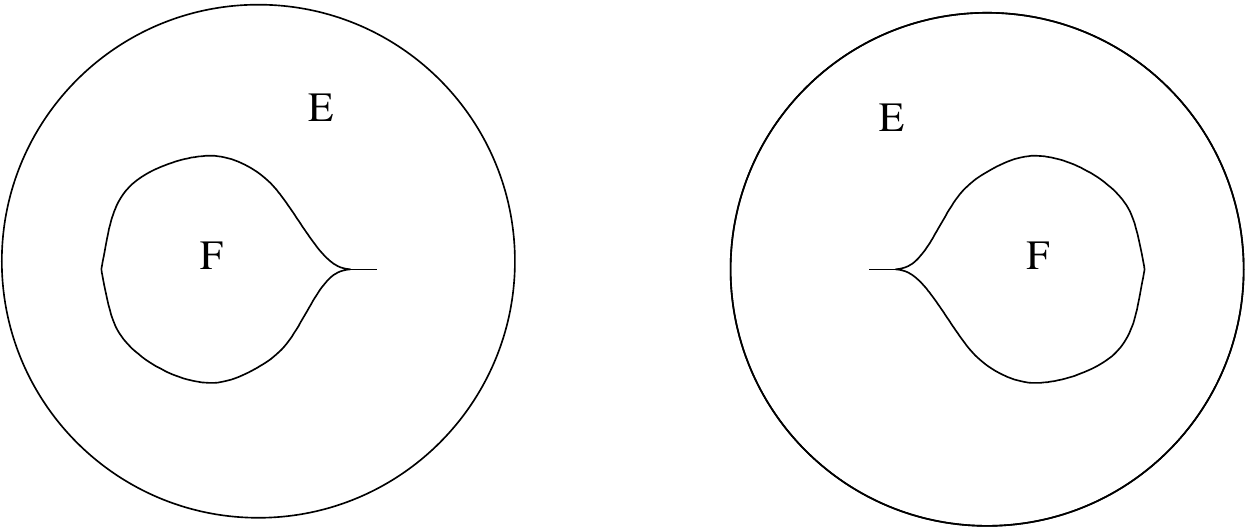}
\caption{Level sets of $u$}
\label{wall1}
\end{center}
\end{figure} 

\begin{figure}[!h]
\begin{center}
\includegraphics[height=14cm]{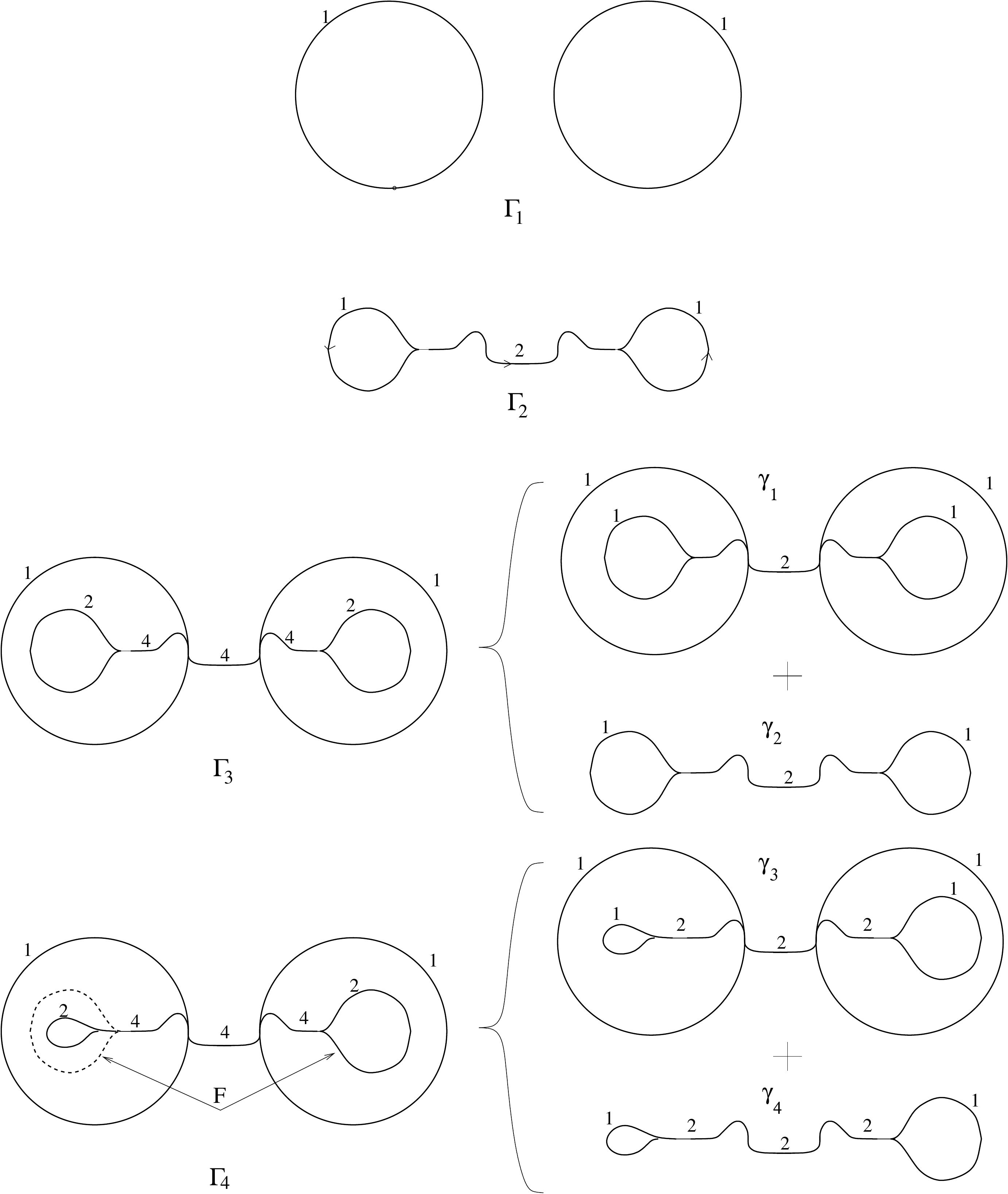}
\caption{The systems $\Gamma_1, \Gamma_2, \Gamma_3, \Gamma_4$ with their multiplicities}
\label{wall2}
\end{center}
\end{figure}
Let $\Gamma_1, \Gamma_2, \Gamma_3, \Gamma_4$ be the systems of curves drawn in Figure~\ref{wall2} together with their multiplicities and consider the piecewise constant function 
$$\Phi(t)=\left\{
\begin{array}{ll}
  \Gamma_1  &  \mbox{if\;} t\in[0,1]\\
  \Gamma_2  &  \mbox{if\;} t\in(1,2]\\
  \emptyset &\text{otherwise}.
\end{array} \right.$$
Clearly, $\Phi\notin  \mathscr{A}(u)$. $\Gamma_1$ and $\Gamma_2$ satisfy Conditions $(i)$, $(ii)$ of Definition \ref{A} but, since the system $\Gamma_1$ does not contain the line joining the two cusp points of $F$, we have $\mathcal{H}^1((\Gamma_2)\setminus \overline{E})\neq 0$ but $\mathcal{H}^1([(\Gamma_2)\setminus\overline{E}]\setminus (\Gamma_1) )\neq 0 $ so  $\Phi$ does not satisfy Condition $(iii)$ of Definition \ref{A}.

\par Consider now the function
$$\Phi(t)=\left\{
\begin{array}{ll}
  \Gamma_3   &\mbox{if\;} t\in[0,1]\\
  \Gamma_2    &\mbox{if\;} t\in]1,2]\\
  \emptyset &\text{otherwise}.
\end{array} \right.$$
where $\Gamma_3$ is built from the curves $\gamma_1$ and $\gamma_2$ together with the multiplicities indicated in Figure~\ref{wall2}. It is easy to check that $\Phi \in \mathscr{A}(u)$ and it must be emphasized that the choice of  $\Phi(t)=\Gamma_2$ for every $t\in ]1,2]$ yields strong geometric constraints. In particular, since the curve joining the two cusp points of $F$ goes out of the set $E$, condition $(iii)$ of Definition~\ref{A} imposes that the trace of $\Phi(t)$ for almost every $t\in [0,1]$ contains $(\Gamma_2)\setminus \overline{E}$. 

\par Let us finally examine the function
$$\Phi(t)=\left\{
\begin{array}{ll}
  \Gamma_4    &\mbox{if\;} t\in[0,1]\\
  \Gamma_2    &\mbox{if\;} t\in(1,2]\\
  \emptyset &\text{otherwise}.
\end{array} \right.$$
where $\Gamma_4$ is built from the curves $\gamma_3$, $\gamma_4$. Remark that, up to a Lebesgue-negligible set, $\Int(\Gamma_4)$ coincides with $E\cup F$ because the multiplicity of the inner curve is everywhere even. In this example,  $\Phi$  satisfies Condition $(iii)$ because the curve joining the two cusp points of $F$ belongs to both $\Gamma_2$ and $\Gamma_4$. However, $\Phi$ does not satisfy Condition $(ii)$ of Definition \ref{A} since $\gamma_4\cap \Int(\Gamma_2)\neq\emptyset$ (see Remark \ref{remarkiii}), so $\Phi\notin  \mathscr{A}(u)$. 
\item {\it Example II}\; Let $v = \mathds{1}_E + \mathds{1}_F$ with $E,F$ like in Figure~\ref{wall3}.
\begin{figure}[!h]
\begin{center}
\includegraphics[height=1.5cm]{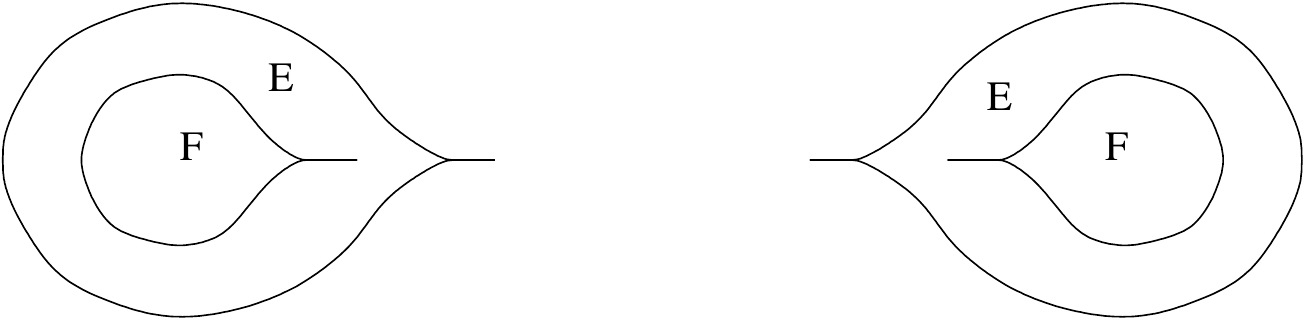}
\caption{Level sets of $v$}
\label{wall3}
\end{center}
\end{figure} 
Let $\Gamma_1, \Gamma_2, \Gamma_3$ be the systems of curves in Figure~\ref{wall4}.
\begin{figure}[!h]
\begin{center}
\includegraphics[height=6cm]{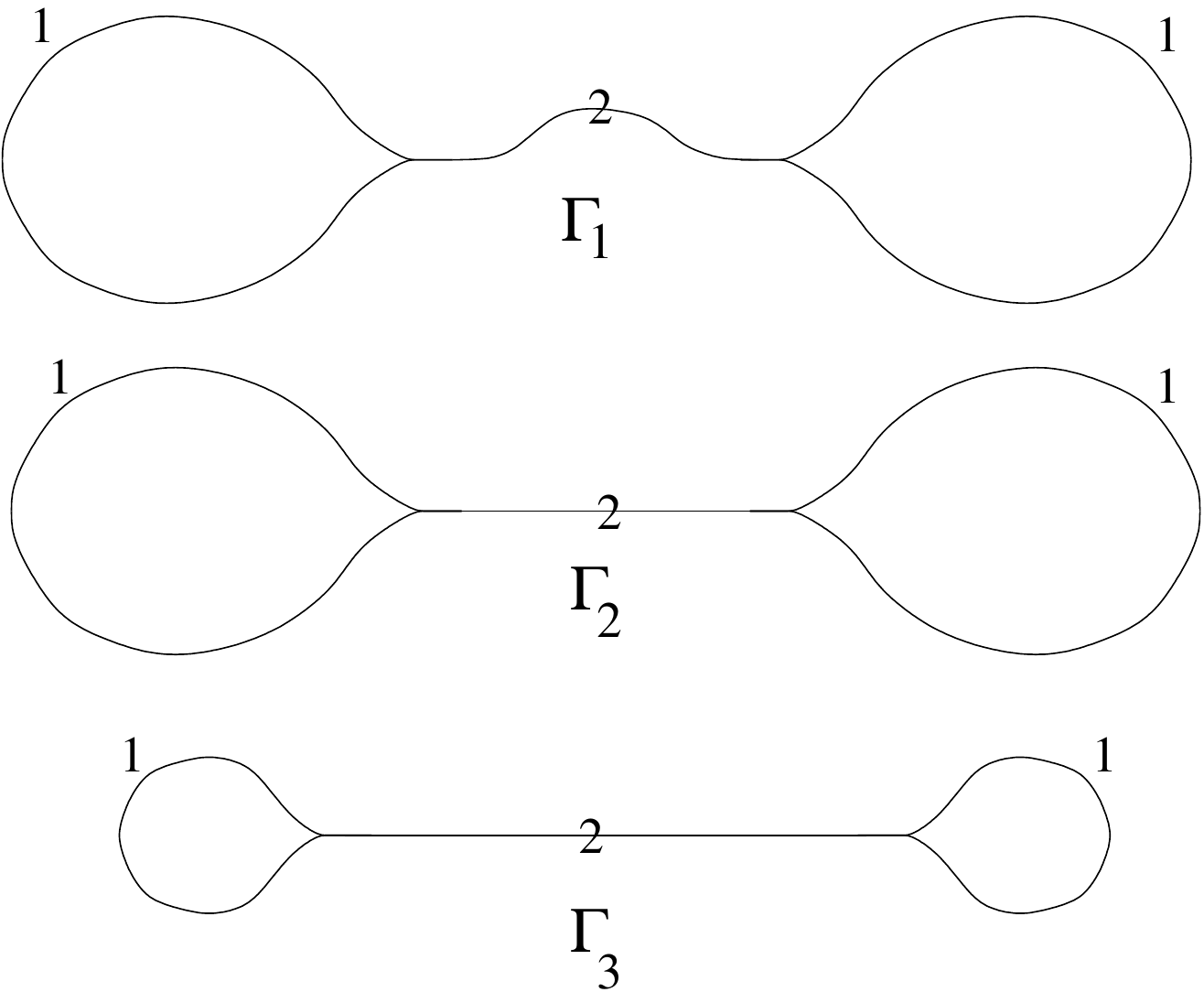}
\caption{The systems $\Gamma_1, \Gamma_2, \Gamma_3$ with their multiplicities}
\label{wall4}
\end{center}
\end{figure} 
\par We consider the following functions: 
$$\Phi(t)=\left\{
\begin{array}{ll}
  \Gamma_1   & \text{if\;} t\in[0,1]\\
  \Gamma_3    &\text{if\;} t\in]1,2]\\
  \emptyset &\text{otherwise}.
\end{array} \right.$$
$$\Psi(t)=\left\{
\begin{array}{ll}
  \Gamma_2    &\text{if\;} t\in[0,1]\\
  \Gamma_3    &\text{if\;} t\in]1,2]\\
  \emptyset &\text{otherwise}.
\end{array} \right.$$

The function $\Phi$ does not satisfy Condition $(iii)$ of Definition~\ref{defiA} because the line joining the two cusps of $F$ goes out of $E$ whereas 
$$\mathcal{H}^1([(\Gamma_2)\setminus \overline{E}]\setminus (\Gamma_1))\neq 0$$
and so  $\Phi \notin \mathscr{A}(u)$. On the contrary, it is easy to check that $\Psi\in \mathscr{A}(u)$ and, again, $\Gamma_2$ imposes strong geometrical constraints for all conditions of Definition~\ref{defiA} to be satisfied.
\end{itemize}
}
\end{oss}

\section{A coarea-type formula for $\oF$}
Recall from the introduction the definition of the functional

$$ G : \begin{array}[t]{lll}
\mathscr{A}&\to &\Bar\R^+\\
\Phi&\mapsto &\ds\int_{\R} W(\Phi(t)) \mbox{d}t.\end{array}$$

Our main result is the representation formula that holds for any $u\in \BV(\R^2)$ with  $\oF(u) < \infty$, that is
$$\oF(u) = \underset{\Phi \in \mathscr{A}(u)} {\Min} G(\Phi).$$

This formula can be easily proved for smooth functions, as shown in the following
\begin{oss}[{\bf The regular case}] 
{\rm Let $u\in \cont^2(\R^2)$ with  $\oF(u) < \infty$. Applying the coarea formula to the system of curves $\Phi[u]$ and since $F$ and $\oF$ coincide for smooth functions, we get immediately that 
$$ \oF(u) = F(u) = G(\Phi[u]).$$
By Definition~\ref{defiA}, $\Phi[u]\in \mathscr{A}(u)$ and
$$G(\Phi[u])\leq G(\Phi) \quad \forall \Phi \in \mathscr{A}(u).$$
Therefore
$$G(\Phi[u])= \underset{\mathscr{A}(u)}{\Min}\; G= \oF(u) =F(u).$$
}
\end{oss}

In order to extend this result for a general function $u$ in $\BV(\R^2)$, let us first address the existence of minimizers of $G$ in $\mathscr{A}(u)$.

\subsection{Existence of minimizers of $G$}\label{G}

The next proposition gives a sufficient condition of compactness with respect to the $\mathscr{A}$-convergence:

\begin{prop}\label{compA}
Let $(\Phi_h)$ be a sequence in $\mathscr{A}$ such that 
$$\underset{h}{\sup} \; G(\Phi_h) < \infty.$$
Then, possibly extracting a subsequence, there exists a function $\Phi \in
\mathscr{A}$ such that 
$$\Phi_h \overset{\mathscr{A}}{\longrightarrow}\Phi \quad \mbox{and} \quad G(\Phi)\leq
\underset{h\rightarrow \infty}{\liminf} \; G(\Phi_h).$$
\end{prop}

\dimo  The proof is essentially the same as the proof of Theorem 2 in \cite{MM} so a few details will be omitted.\\
{\bf Step 1 : Convergence of the energies $W(\Phi_h(t))$.}

Let $N \in \mathbb{N}, k\in \mathbb{Z}$ and let us consider the dyadic intervals on
$\R$ :
$$I_{N,k}=[k2^{-N}, (k+1)2^{-N}[.$$
We define the functions
$$f_h^N(t)=2^N\int_{I_{N,k}} W(\Phi_h(s)) \mbox{d}s $$
where $I_{N,k}$ is the unique dyadic interval containing $t$. The function
$f_h^N(t)$ is constant on each interval $I_{N,k}$ and for every $t\in \R$ we
have $f_h^N(t)\leq 2^N G(\Phi_h)$. So, by a diagonal extraction, we can take
a subsequence (not relabelled) such that 
$$\forall (N,k) \quad f_h^N(t)\rightarrow f^N(t) \quad \forall t\in I_{N,k}.$$ 
Moreover we can write 
$$I_{N,k}=[k2^{-N}, (2k+1)2^{-N-1}[\cup [(2k+1)2^{-N-1}, (k+1)2^{-N}[$$
and 
$$f_h^{N+1}(t) = 2^{N+1}\int_{I_{N+1,2k}}W(\Phi_h(s)) \mbox{d}s \quad \forall t\in [k2^{-N}, (2k+1)2^{-N-1}[$$
$$f_h^{N+1}(t) = 2^{N+1}\int_{I_{N+1,2k+1}}W(\Phi_h(s)) \mbox{d}s \quad \forall t\in [(2k+1)2^{-N-1}, (k+1)2^{-N}[,$$
therefore
$$\int_{I_{N,k}}f_h^{N+1}(s) \mbox{d}s =  \int_{I_{N,k}}W(\Phi_h(s)) \mbox{d}s$$
$$f_h^N(t) = 2^N \int_{I_{N,k}}f_h^{N+1}(s)\mbox{d}s \quad \forall t\in I_{N,k}.$$
Then, by the Dominated Convergence Theorem, we get
\begin{equation}\label{marti}
 f ^N(t) = 2^N \int_{I_{N,k}}f^{N+1}(s) \mbox{d}s
\end{equation}
and in addition, by Fatou's Lemma, 

\begin{equation}\label{Doob}
 \int_{\R}f^N(s) \mbox{d}s \leq \underset{h\rightarrow \infty}{\liminf} \;\; \int_{\R}f_h^N(s)\mbox{d}s\leq  \underset{h}{\sup} \; G(\Phi_h).
\end{equation}

\eqref{marti} and \eqref{Doob} show that $(f^N)$ is a bounded positive martingale thus, by the convergence theorem for martingales (\cite[Thm 2.2, p. 60]{RY}, there exists $f\in \lp{1}(\R)$ such that $f^N \rightarrow f$ a.e.

\par {\bf {Step 2 : Definition of a limit system of curves $\Phi$.}}

Let $N\in \mathbb{N}, k\in \mathbb{Z}$. We have
\begin{equation}\label{lemmino}
\underset{I_{N,k}}{\sup}\; \underset{h}{\sup}\; f_h^N < \infty.
\end{equation}

\begin{lemma}\cite{MM}\label{lemmatec}
Let $A_h:= \left\lbrace t\in I_{N,k} : W(\Phi_h(t)) \leq 2^N\int_{I_{N,k}}
W(\Phi_h(s)) \mbox{d}s + \dfrac{1}{N}\right\rbrace $. Then there exists $t_{N,k}\in
I_{N,k}$ such that, possibly passing to a subsequence, 
$$t_{N,k}\in A_h, \quad\quad \forall h\in\mathbb{N}.$$
\end{lemma}

\dimo see~\cite{MM} Lemmas 4 and 5
\qed

For every dyadic interval $I_{N,k}$ we consider the real number $t_{N,k}$  given by
the previous lemma then, possibly extracting a subsequence, we have 
$$W(\Phi_h(t))\leq 2^N\int_{I_{N,k}} W(\Phi_h(s)) \mbox{d}s + \dfrac{1}{N}\quad
\forall h$$
and so by the compactness theorem in $\wpq{2}{p}$ there exists a subsequence (not
relabelled) and a limit system of curves $\Phi(t_{N,k})$ such that 
$$\Phi_h(t_{N,k}) \overset{\wpq{2}{p}}{\rightharpoonup}  \Phi(t_{N,k})\quad
\mbox{and} \quad W(\Phi(t_{N,k}))\leq \underset{h\rightarrow
\infty}{\liminf}\; W(\Phi_h(t_{N,k})).$$
Remark that, since the curves of $\Phi_h(t_{N,k})$ are without crossing, because of the $\cont^1$-convergence the curves of $\Phi(t_{N,k})$ are without crossing as well. 
\par Since the $I_{N,k}$'s are countably many, we can use a diagonal extraction argument to find a subsequence, still denoted by $\Phi_h$, and a limit system of curves of class $\wpq{2}{p}$, denoted by $\Phi(t_{N,k})$, such that for each $t_{N,k}$ given by
 the previous lemma:
$$\Phi_h(t_{N,k}) \overset{\wpq{2}{p}}{\rightharpoonup}  \Phi(t_{N,k})\quad
\mbox{and}\quad  W(\Phi(t_{N,k}))\leq \underset{h\rightarrow
\infty}{\liminf}\; W(\Phi_h(t_{N,k})).$$
\par Furthemore for every $t_{N,k}> t_{N, k'}$ the systems $\Phi_h(t_{N,k})$ and $\Phi_h(t_{N,k'})$  are without crossing and $\Int(\Phi_h(t_{N,k}))\subseteq \Int(\Phi_h(t_{N,k'}))$  so, because of the $\cont^1$ convergence, also  $\Phi(t_{N,k})$ and $\Phi(t_{N,k'})$ are without crossing and $\Int(\Phi(t_{N,k}))\subseteq \Int(\Phi(t_{N,k'}))$.\\

Let us now  see how a limit curve  can be defined for every $t$. Let $t\in
\R, N\in \mathbb{N}$ so there exists  $k_N \in\mathbb{Z}$ such that $t\in
I_{N, k_N}$. We have 
$$W(\Phi(t_{N,k_N}))\leq f^N(t) + \dfrac{1}{N}\overset{N \rightarrow
\infty}{\rightarrow} f(t)\quad \mbox{a.e.}.$$
Then by the weak compactness of $\wpq{2}{p}$ there exists a subsequence, still denoted by $(\Phi(t_{N,k_N}))$, and a system of curves of class $\wpq{2}{p}$, denoted by $\Phi(t)$, such that 
$$\Phi(t_{N,k_N})
\overset{\wpq{2}{p}}{\rightharpoonup} \Phi(t)\quad \mbox{and} \quad
W(\Phi(t))\leq \underset{h\rightarrow \infty}{\liminf}\; W(\Phi(t_{N,k_N}))\leq
f(t).$$

This procedure can be applied for  almost every $t\in \R$ so that, if we can prove that $\Phi\in \mathscr{A}$, we will conclude that $\Phi_h \overset{\mathscr{A}}{\rightarrow} \Phi$. Observe that
\begin{itemize}
\item  the curves of $\Phi(t_{N,k})$ are without crossing, as was shown before, and 
$${\mbox{Int}}(\Phi(t_{N,k}))\subseteq {\mbox{Int}}(\Phi(t_{N,k'}))$$ for every $t_{N,k}> t_{N, k'}$;

\item for almost every $\underline{t}, \overline{t}\in \mathbb{R}$, $\underline{t}\neq\overline{t} $,  we can find $\overline{t}_{N,k_N} \rightarrow \overline{t}$ and 
$\underline{t}_{N,k_N'} \rightarrow \underline{t}$ such that 
$$\Phi(\overline{t}_{N,k_N})
\overset{W^{2,p}}{\rightharpoonup}  \Phi(\overline{t}),\;\; \Phi(\underline{t}_{N,k_N'})
\overset{W^{2,p}}{\rightharpoonup}  \Phi(\underline{t}).$$ Since $\Phi(\overline{t}_{N,k_N})$ and $ \Phi(\underline{t}_{N,k_N'})$ are without crossing, because of the $\cont^1$ convergence we get that $\Phi(\overline{t})$ and $\Phi(\underline{t})$ are without crossing. Moreover we can suppose, for $N$ large enough, $\overline{t}_{N,k_N} > \underline{t}_{N,k_N'}$ and since $ {\mbox{Int}}(\Phi(\overline{t}_{N,k_N})) \subseteq {\mbox{Int}}(\Phi(\underline{t}_{N,k_N'})) $, the  $\cont^1$ convergence implies that
$$ {\mbox{Int}}(\Phi(\overline{t})) \subseteq {\mbox{Int}}(\Phi(\underline{t}));$$

\item we have to prove that  for almost every  $\underline{t}, \overline{t}\in \R$, $\underline{t}\neq\overline{t} $ the system $\Phi(t)$ satisfies condition~$(iii)$ of Definition~\ref{A}.
\par For every $t$ we let
$$E(\Phi(t)) = \R^2\setminus\overline{\Int(\Phi(\underline{t}))} .$$
By contradiction, suppose that there exists $\overline{t}\geq\underline{t}$ and $\gamma_{\overline{t}}^i  \in \Phi(\overline{t})$ such that 
$$\mathcal{H}^1\left(  (\gamma_{\overline{t}}^i)  \cap E(\Phi(\underline{t})) \right)\neq 0$$
and
$$\mathcal{H}^1\left(  [(\gamma_{\overline{t}}^i)  \cap E(\Phi(\underline{t}))] \setminus(\Phi(\underline{t})) \right) \neq 0.$$
Then we can find two sequences $\{(h_N, k_N)\}, \{(h_N, k_N')\}$ such that $t_{N,k_N} \rightarrow \overline{t}$ and 
$t_{N,k_N'} \rightarrow \underline{t}$, with $t_{N,k_N'} < t_{N,k_N}$ for every $N$ large enough, that satisfy
\begin{equation}\label{condiii}
 \Phi_{h_N}(\underline{t}_{k_N,N})\overset{\cont^1}{\longrightarrow}
\Phi(\underline{t})\quad \mbox{and} \quad
 \gamma_{h_N}\overset{\cont^1}{\longrightarrow} \gamma_{\overline{t}}^i\;, \quad \gamma_{h_N} \in \Phi_{h_N}(\overline{t}_{k_N,N}).
\end{equation}

Because of the $\cont^1$-convergence, for $N$ large enough we have 
$$\mathcal{H}^1\left(  (\gamma_{h_N})  \cap E(\Phi(t_{N,k_N'})) \right)\neq 0$$
and
$$\mathcal{H}^1\left(  [(\gamma_{h_N}) \cap E(\Phi(t_{N,k_N'}))] \setminus (\Phi(t_{N,k_N'})) \right) \neq 0$$
which gives a contradiction with the fact that condition $(iii)$ holds for the functions $\Phi_{h_N}$.
\end{itemize}

\par Finally, we have defined a collection of curves $\Phi \in \mathscr{A}$ such that $\Phi_h
\overset{\mathscr{A}}{\longrightarrow}\Phi$. Moreover, 
$$G(\Phi)=\int_{\R} W(\Phi(t)) \mbox{d}t \leq \int_{\R} f(t)
\mbox{d}t \leq \underset{N\rightarrow \infty}{\liminf}\int_{\R} f^N(t)
\mbox{d}t$$
$$\leq  \underset{h\rightarrow \infty}{\liminf}\;\underset{N\rightarrow
\infty}{\liminf}\int_{\R} f_h^N(t) \mbox{d}t=  \underset{h\rightarrow
\infty}{\liminf}\int_{\R} W(\Phi_h(t)) \mbox{d}t =\underset{h\rightarrow
\infty}{\liminf}\; G(\Phi_h).$$
\qed

The next theorem states the existence of minimizers to $G$ in $\mathscr{A}(u)$.

\begin{thm}\label{minG}
Let $u\in \BV(\R^2)$ with $\oF(u) < \infty$. Then $\mathscr{A}(u)\neq \emptyset$, the problem
$$\underset{\Phi \in \mathscr{A}(u)}{\Min} G(\Phi) $$
has a solution, and
$$\underset{\Phi \in \mathscr{A}(u)}{\Min} G(\Phi)  \leq \oF(u).$$
\end{thm}

\dimo Let $\{u_h\} \subset \cont^2(\R^2)$ be a sequence converging to $u$ in
$\lp{1}$ such that $F(u_h)$ is uniformly bounded and $\oF(u)=\underset{h\rightarrow
\infty}{\lim}F(u_h)$. 
We can associate to every $u_h$ the function $\Phi_h=\Phi[u_h]$
and by the coarea formula  we have
$$F(u_h)=G(\Phi_h),\quad \underset{h}{\sup}\;G(\Phi_h)<\infty.$$
Then, by Proposition \ref{compA}, there exists a subsequence (not relabelled) and $\Phi\in \mathscr{A}$
such that 
$$\Phi_h \overset{\mathscr{A}}{\longrightarrow} \Phi ,\quad \quad G(\Phi) \leq \underset{h \rightarrow \infty}{\liminf}\; G(\Phi_h).$$
Let us now prove that $\Phi\in\mathscr{A}(u)$. By definition of the $\mathscr{A}$-convergence we have
$$ |\Int(\Phi_{h}(t_{N, k_N}))\Delta  \Int(\Phi(t_{N, k_N}))| = \arrowvert \{ u_{h} > t_{N, k_N}\}\Delta  \Int(\Phi(t_{N, k_N}))\arrowvert \rightarrow 0\quad \mbox{if}\quad
h \rightarrow \infty$$
and, since $\arrowvert \{ u_{h} > t\} \Delta  \{u>t\}\arrowvert \rightarrow 0$ for $\mathcal{L}^1$-almost every $t$, in Lemma \ref{lemmatec} we can choose $t_{N, k_N}$ such that $\arrowvert \{ u_{h} > t_{N, k_N}\} \Delta  \{u>t_{N, k_N}\}\arrowvert \rightarrow 0$ and we have 
\begin{equation}\label{3}
 \{u>t_{N, k_N}\}=\Int(\Phi(t_{N, k_N})) \quad \mbox{(up to a Lebesgue negligible  set)} \quad \forall N.
\end{equation}
In addition, for almost every $t$, $\Phi(t)$ is the
weak $\wpq{2}{p}$ limit  of a sequence  $\{\Phi(t_{N, k_N})\}$ where $t_{N,k_N}\rightarrow
t$ as $N\rightarrow \infty$, and it follows that
$$|\Int(\Phi(t_{N, k_N}))\Delta  \Int(\Phi(t))|\rightarrow 0\quad \mbox{if}\quad
N \rightarrow \infty.$$
Now, the interiors of the systems $\Phi(t)$ are nested, and so, using  Lemma~\ref{mon} for the family $(\{u>t\})_t$ and~\eqref{3}, we get 
\begin{equation}\label{2}
 \{u>t\}=\Int(\Phi(t)) \quad \mbox{(up to a Lebesgue negligible  set)} \quad \text{for a.e. $t$}.\end{equation}
Being $u$ of bounded variation, it follows that for almost every $t$
$$\partial^* \{u>t\} =\partial^*\Int(\Phi(t)) \quad\text{up to a $\mathcal{H}^1$-negligible  set}$$
therefore
$$\partial^* \{u>t\} \subset\Phi(t)\quad\text{up to a $\mathcal{H}^1$-negligible  set.}$$
This proves that there exists $\Phi\in\mathscr{A}(u)$ such that $G(\Phi)\leq \oF(u)$.

To show the existence of minimizers to $G$ in $\mathscr{A}(u)$, it suffices to take $(\Phi_h)\subset \mathscr{A}(u)$ a minimizing sequence, i.e.
$$\underset{\mathscr{A}(u)}{\inf}\, G=\underset{h\rightarrow \infty}{\lim}G(\Phi_h)$$
and  $\underset{h}{\sup}\;G(\Phi_h)< \infty$. Using exactly the same argument as above, we can conclude that there exists $\Phi\in \mathscr{A}(u)$ such that 
$$\Phi_h \overset{\mathscr{A}}{\longrightarrow} \Phi ,\quad \quad G(\Phi) \leq \underset{h \rightarrow \infty}{\liminf}\; G(\Phi_h)=\underset{\mathscr{A}(u)}{\inf}\, G$$
therefore 
$$\underset{\Phi \in \mathscr{A}(u)}{\Min} G(\Phi) $$
has a solution, and
$$\underset{\Phi \in \mathscr{A}(u)}{\Min} G(\Phi)  \leq \oF(u).$$
\qed

\subsection{Connections between $G$ and $\oF$}\label{formula}

We can deduce from the previous theorem a first representation result for $\oF$:

\begin{lemma}\label{regolare}
 Let $E \subset  \R^2$ be a bounded open set with $\partial E \in  \wpq{2}{p}$ and let $u=c\mathds{1}_E$, $c >0$. Then 
$$\oF(u) =  c \oW( E)=c\int_{\partial E}(1+|\kappa_{\partial E}|^p)d\hone.$$
\end{lemma}

\dimo By Corollary~3.2 in~\cite{BDP}, $c\oW(E)=G(\Psi)=c\ds\int_{\partial E}(1+|\kappa_{\partial E}|^p)d\hone$ where $\Psi$ is the map defined by $\Psi(t) = \partial E$ for every $t\in[0,c]$, $\emptyset$ otherwise. Therefore, if $\oF(u)<+\infty$ then, by Theorem~\ref{minG} and Definition~\ref{defiA},
$$\oF(u) \geq  \underset{\mathscr{A}(u)}{\Min}\,G = G(\Psi) = c \oW( E).$$
The fact that $\oF(u)<+\infty$ and the reverse inequality $\oF(u) \leq c \oW( E)$ will follow if we can find a sequence $\{u_h\}\subseteq \cont^2(\R^2)$ such that $u_h \rightarrow u$ in $\lp{1}$ and 
$$F(u_h)=G(\Phi[u_h]) \rightarrow G(\Psi).$$
Since $c\oW(E)=G(\Psi)<+\infty$, there exists a sequence of open sets $(E_m)$ of class $\cont^2$ such that 
\begin{equation}
|E_m\Delta E|\to 0\qquad\text{and}\qquad cW(E_m)\to G(\Psi)\label{approxBDP-G}
\end{equation}
\par
\begin{figure}[!h]
\begin{center}
\includegraphics[height=6cm,angle=270]{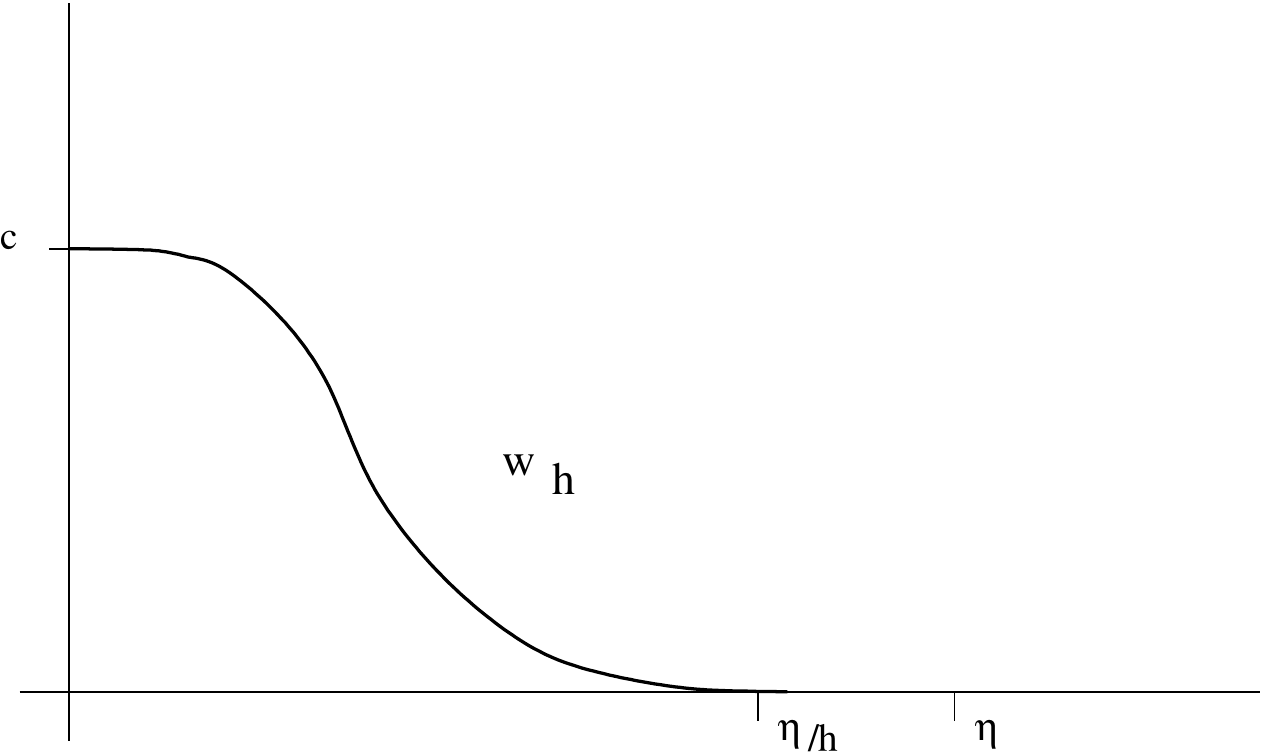}
\caption{Graph of the cut-off function $w^m_h$}
\label{eta}
\end{center}
\end{figure}

 Let $m\in\N$. From the properties of the distance function (see \cite[\S 14.6, p .354]{GT}), we can find $\eta> 0$ such that $ \,d(x):={\rm{dist}}(x, \partial E_m) \in \cont^2(E^\eta_m)$ where $E^\eta_m= \{x\in \R^2 \setminus \overline{E_m} : \,d(x)< \eta \}$ and the curvature of $\partial E_m$ is bounded by $\eta^{-1}$ so that 
\begin{equation}\label{curvbound}
 1+\,d(x) k_{\partial E_m}(x) \neq 0 \quad \quad  \forall x \in E^\eta_m.
\end{equation}
For every $h\in\N^*$ we consider the cut-off function $w^m_h : [0, \eta/h] \rightarrow [0,c]$, $w^m_h\in \cont^2([0, \eta/m])$ whose graph is represented in Figure~\ref{eta}. 
Then we define the sequence of smooth functions:
$$u^m_h(x)=\left\{
\begin{array}{ll}
   c         & \text{if } x\in \overline{E_m},\\
   w^m_h(\,d(x))   & \text{if } x\in E^{ \eta/h}_m, \\
   0          & \text{otherwise.}
\end{array} \right.$$
Thus $u^m_h \overset{\lp{1}}{\rightarrow} u^m:=\one{E_m}$, as $h \rightarrow +\infty$, and for every $t\in [0,c]$ we get
\begin{equation}\label{livdist}
 \Gamma_t^{m,h} = \{x \in \R^2 : u^m_h(x) = t\}=\left\lbrace x \in \R^2 :  \,d(x)= (w^m_h)^{-1}(t) \right\rbrace.
\end{equation}
Therefore, for $h$ large enough, $\Gamma_t^{m,h}$ can be parametrized as (see~\cite[\S 5.7, p.115]{G}):
$$\Gamma_t^{m,h}(y)= y + \delta n(y), \quad y\in \partial E_m ,\quad \delta= (w^m_h)^{-1}(t)$$
where $n(y)$ is the outer unit normal to $\partial E_m$ at $y$. Using a positively-oriented arc-length parameterization $s\mapsto\alpha(s)$ of $\partial E_m$, we can parametrize $\Gamma_t^{m,h}$ as (the dependence on $t$ is omitted to simplify the notations)
$$\gamma(s) = \alpha(s) + \delta n(\alpha(s))= \alpha(s) - \delta J \alpha'(s)$$
where $J$ is the rotation operator $J(x,y)=(-y,x)$, $(x,y) \in \R^2.$
Thus
$$\gamma'(s)= \alpha'(s) -\delta J\alpha''(s)= \alpha'(s) - \delta k_\alpha(s) J^2\alpha'(s) = \left(1+\delta k_\alpha(s)\right)\alpha'(s),$$
where, with a small abuse of notation, $k_\alpha(s)$ now denotes the signed scalar curvature instead of the vector curvature, i.e. $\alpha''(s)=k_\alpha(s)\vec{N}$ being $(\alpha'(s),N)$ a direct frame with $|\vec{N}|=1$.
Hence $\gamma$ is regular at any $s$ such that $1+\delta k_\alpha(s) \neq 0$ and 
$$ \gamma''(s)= \left(1+\delta k_\alpha(s)\right) k_\alpha(s) J\alpha'(s) +\delta k_\alpha'(s)\alpha'(s)$$
and so
$$k_\gamma(s)= \dfrac{\langle \gamma'', J \gamma' \rangle}{\|\gamma' \|^3}= \dfrac{k_\alpha(s)}{|1+\delta k_\alpha(s)|}.$$
Then, by \eqref{curvbound}, $\Gamma_t^{m,h}$ is smooth for $h$ large enough and 
$$k_{\Gamma_t^{m,h}}(\Gamma_t^{m,h}(y))= \dfrac{k_{\partial E_m}}{|1+\delta k_{\partial E_m}|}(y)$$
$$\mbox{d}\mathcal{H}^1 \res \Gamma_t^{m,h} = |1+\delta k_{\partial E_m}| \mbox{d}\mathcal{H}^1 \res \partial E_m $$
and we  get
$$ W(\Gamma_t^{m,h})= \int_{ \Gamma_t^{m,h}}[1+|k_{\Gamma_t^{m,h}}|^p] \mbox{d}\mathcal{H}^1 =\int_{ \partial E_m}\left[1+ \left|\dfrac{k_{\partial E_m}}{1+\delta k_{\partial E_m}}\right|^p \right]  |1+\delta k_{\partial E_m}| \mbox{d}\mathcal{H}^1.$$
Then, remark that
\begin{equation}\label{lunghW}
 |\mathcal{H}^1(\Gamma_t^{m,h}) - \mathcal{H}^1(\partial E_m)|\leq \int_{\partial E_m} \delta |k_{\partial E_m}| \mbox{d}\mathcal{H}^1  \leq C_1\delta W( E_m).
\end{equation}
Moreover,
$$\left|\int_{\Gamma_t^{m,h}} |k_{\Gamma_t^{m,h}}|^p \mbox{d}\mathcal{H}^1 - \int_{\partial E} |k_{\partial E_m}|^p \mbox{d}\mathcal{H}^1\right| \leq \int_{\partial E_m}|k_{\partial E_m}|^p\left| \dfrac{1}{|1+\delta k_{\partial E_m}|^{p-1}} -1\right|  \mbox{d}\mathcal{H}^1$$
thus, by Lebesgue's Dominated Convergence Theorem, 
\begin{equation}\label{curvW}
 \left|\int_{\Gamma_t^{m,h}} k_{\Gamma_t^{m,h}}^p \mbox{d}\mathcal{H}^1 - \int_{\partial E_m} k_{\partial E_m}^p \mbox{d}\mathcal{H}^1\right|\rightarrow 0 \quad \quad \text{as} \quad \delta\rightarrow 0 \quad \quad (\text{i.e.}\quad  h\rightarrow \infty).
\end{equation}
It follows from \eqref{lunghW} and \eqref{curvW} that
\begin{equation}\label{distW}
 \underset{h \rightarrow \infty}{\lim} |W(\Gamma_t^{m,h}) - W( E_m)| =0 \quad \forall t\in [0,c].
\end{equation}
Consider now a subsequence $\Phi_k=\Phi[u^{m(k)}_{h(k)}]$. We have $\Phi_k(t)= \Gamma_t^{m(k),h(k)}$ for every $k$ and every $t\in [0,c]$. In addition, 
$$|F(u^{m(k)}_{h(k)}) - G(\Psi)|= |G(\Phi[u^{m(k)}_{h(k)}]) - G(\Psi)|\leq   \int_0^{c} \left| W(\Gamma_t^{m(k),h(k)})-\oW( E) \right|dt $$
Using~\eqref{approxBDP-G},~\eqref{distW} and a diagonal extraction argument, we can find a subsequence such that, keeping the same labeling and applying the Dominated Convergence Theorem
$$ u^{m(k)}_{h(k)} \overset{\lp{1}}{\longrightarrow} u \quad \text{and} \quad F(u^{m(k)}_{h(k)}) \rightarrow G(\Psi) = c\oW( E).$$
Since $\oF(u)\leq\liminf_{k\to\infty}F(u^{m(k)}_{h(k)})$ the conclusion follows.
\qed

This proof illustrates that, at least in the case of $u$ being the characteristic function of a $\wpq{2}{p}$ set, the equivalence between $\oF$ and $G$ follows from a smoothing argument together with a control of the energy. The purpose of the next lemma is to show that a similar strategy also holds for more complicated  functions. This lemma is essentially the same as Lemma 6 in~\cite{MM} so most details will be omitted.

\begin{lemma}\label{approx}
Let $u\in \BV(\R^2)$ with $\oF(u)<\infty$ and let $\Phi \in
\mathscr{A}(u)$. Then, for every
$\eta> 0$, there exists  $\tilde{u}\in \BV(\R^2)$  such that
$$\|u-\tilde{u}\|_{\lp{1}} \leq \eta,$$
$\partial^* \{ \tilde{u} > t \}$ is, for almost every $t$, a finite system of curves of class $\wpq{2}{p}$ without contact or auto-contacts, and any two $\partial^* \{ \tilde{u} > t \}$ and $\partial^* \{ \tilde{u} > s \}$ are disjoint (pointwisely) for almost every $s\not=t$. In addition, the function $\tilde{\Phi}$ defined as
$$t\mapsto \tilde{\Phi}(t)= \partial^* \{ \tilde{u} > t \} $$
belongs to $\mathscr{A}(\tilde{u})$ and 
$$\arrowvert G(\Phi)-G(\tilde{\Phi})\arrowvert \leq \eta .$$
\end{lemma}

\dimo  We suppose $G(\Phi)\neq 0$. The idea of the proof  is to move  smoothly every system of curves $\Phi(t)$ to get a new family of systems of curves of class $\wpq{2}{p}$ belonging to $\mathscr{A}$ with no contact or auto-contact between curves, and with an energy close to the energy of $\Phi$. Then a function of bounded variation can be canonically defined.


\par Let $\{t_n\}_{n\geq 1}$ denote a countable and dense subset of $\R$ such that $\Phi(t_n)$ is well-defined for every $n\geq 1$. Following the proof of Lemma 6 in~\cite{MM} one associates with each finite system of curves $\Phi(t_n)$ a smooth operator $\T_n$ that separates every curve of  $\Phi(t_n)$ from all other curves of $\Phi$ and that removes the auto-contacts by separating the corresponding arcs. $\T_n$ is chosen so that the energy of all curves that have been moved does not increase too much. It is shown in~\cite{MM} that the limit operator $\T=\lim_{n\to\infty}\T_n\circ \T_{n-1}\circ\cdots\circ \T_1$ is well-defined and smooth,  that all curves of $\T(\Phi)$ are without contact or auto-contact, and that, given any $\eta>0$, $\T$ can be designed so that
\begin{equation}\label{stimaii}
|G(\T(\Phi))- G(\Phi)|\leq \eta
\end{equation}

Furthermore, it is easy to check that the separation process preserves the conditions of Definition~\ref{A} thus $\tilde{\Phi}:=\T(\Phi)\in\mathscr{A}$. In particular
\begin{equation}\label{nested}
  \Int(\tilde{\Phi}(s)) \subseteq \Int(\tilde{\Phi}(t))\quad\text{whenever $s>t$}.
\end{equation}
By standard arguments, the function defined by $\tilde{u}(x)=\sup \{t : x\in \Int(\tilde{\Phi}(t)) \}$ is measurable and
$$ \Int(\tilde{\Phi}(s)) \subseteq  \{\tilde{u}> t\} \subseteq \Int(\tilde{\Phi}(\ell))\qquad\text{when $\ell<t<s$}.$$ 
Let $F\subseteq \R$ be the countable set of discontinuities, with respect to the Lebesgue measure, of the monotone family $\{\Int(\tilde{\Phi}(t))\}$. Taking  $s_k \searrow t$ and $l_k \nearrow t $, Lemma \ref{mon} implies that 
$$|\{\tilde{u}>t\}\Delta \Int(\tilde{\Phi}(t))|=0 \quad \forall t \in \R\setminus F.$$

\par Let us now prove that $\tilde{u}\in \BV(\R^2)$. Given $\eta>0$ and using the convention $\T_0=\text{Id}$, we can redefine the separation operator at step $i\geq 1$ so that, by the coarea formula,
$$\arrowvert \Int(\T_i\circ\cdots\circ\T_1(\Phi)(t))  \Delta \Int(\T_{i-1}\circ\cdots\circ\T_0(\Phi)(t)) \arrowvert < \varepsilon_i = \eta 2^{-i}\hone(\Phi(t))$$
for almost every $t$, therefore
$$\arrowvert  \Int(\tilde{\Phi}(t)) \Delta \{u>t\}\arrowvert = \arrowvert   \{\tilde{u}>t\}\Delta \{u>t\}\arrowvert  <  \eta\hone(\Phi(t)).$$
It follows that 
$$\|u-\tilde{u}\|_{\lp{1}} = \int_{-\infty}^{+\infty}\int_{\R^2}\arrowvert\mathds{1}_{\{u>t\}}-\mathds{1}_{\{\tilde{u}>t\}}\arrowvert\, dx\, dt < \eta\,\int_{-\infty}^{+\infty}\hone(\Phi(t))\,dt\leq\eta G(\Phi).$$
In particular, $\tilde{u}\in \lp{1}(\R^2)$. Besides, all curves in $\tilde\Phi$ are mutually disjoint and without auto-contacts and   $\tilde{u}$ is continuous thus 
$$\partial^*\{\tilde{u}> t\}= \partial [\Int\tilde{\Phi}(t)] = \tilde{\Phi}(t) \quad \mbox{ (up to a $\mathcal{H}^1$-negligible  set)}.$$
We know from  \eqref{stimaii} that $G(\tilde\Phi)<+\infty$ therefore $t\mapsto \mathcal{H }^1(\partial^*\{\tilde{u}> t\})$ is in $\lp{1}(\R)$ thus, since $\tilde u\in\lp{1}(\R^2)$, $\tilde{u}\in \BV(\R^2)$ by the coarea formula.
\qed

We can now state our main result :

\begin{thm}\label{principale}
 Let $u\in \BV(\R^2)$ with  $\oF(u) < \infty$. Then
$$\oF(u) = \underset{\Phi \in \mathscr{A}(u)} {\Min} G(\Phi).$$
\end{thm}
\dimo
 Writing  $u=u^+-u^-$  where $u^+(x)=\max\{u(x),0\}$ and $u^-(x)=\max\{-u(x),0\}$, and observing that $\oF(u)\leq \oF(u^+)+\oF(u^-)$,  we can suppose $u\geq 0$. In view of Theorem~\ref{minG}, we only have to prove that
\begin{equation}\label{maggiorazione}
 \underset{\Phi \in\mathscr{A}(u)}{\Min} G(\Phi)  \geq \oF(u).
\end{equation}
Let us show that, for every $\varepsilon > 0$, we can find $v\in \cont^2(\R^2)$ such that
$$\|u-v\|_{\lp{1}} \leq \varepsilon$$
$$\arrowvert G(\Phi)-F(v) \arrowvert\leq \varepsilon $$
where $\Phi$ is a minimizer of $G$ on $\mathscr{A}(u)$.
\par For every $\varepsilon >0$, by   Lemma \ref{approx}, there exists $\tilde{u}\in
\BV(\R^2)$, with $\wpq{2}{p}$ level lines without contacts or self-contacts, such that 

\begin{equation}\label{l1}
 \|u-\tilde{u}\|_{\lp{1}} \leq \varepsilon/4
\end{equation}
\begin{equation}\label{w1}
 \arrowvert G(\Phi)- G(\tilde{\Phi})\arrowvert \leq \varepsilon/4
\end{equation}
where $\tilde{\Phi}(t)= \partial^* \{\tilde{u} > t\}$ for almost every $t\in\R$

Moreover we can find a set $\{t_n\}_{n\in\N}$ with $t_0=0$,  such that, defining
$$\tilde{v}(x)= \sum_{n\in\N} (t_{n+1}-t_n)\mathds{1}_{\{\tilde{u}>t_{n}\}}(x),$$
there holds
\begin{equation}\label{l2}
 \|\tilde{u}-\tilde{v}\|_{\lp{1}} \leq \varepsilon/4,
\end{equation}
\begin{equation}\label{w2}
\text{and}\qquad \arrowvert G(\tilde{\Phi})- G(\Phi_{\tilde{v}})\arrowvert \leq \varepsilon/4,
\end{equation}
where $\Phi_{\tilde{v}}$, defined by $\Phi_{\tilde{v}}(t)=\partial^* \{\tilde{v} > t\}$ for almost every $t$, satisfies 
$$G(\Phi_{\tilde{v}})=\sum_{n\in\N}(t_{n+1}-t_n)W(\partial^* \{\tilde{u} > t_{n}\}).$$

Therefore we have approximated $\tilde{u}$ in $\lp{1}$ by a piecewise
constant function $\tilde{v}$ whose level lines are systems of curves of class  $\wpq{2}{p}$ without self-contacts and  with finite
$p$-elastica energy. Remark that 
$$\{\tilde{v}> t_n\}= \{\tilde{u}>t_{n}\}\quad(\operatorname{mod}\;{\cal L}^2)\qquad \forall n\in\N.$$
and
$$\partial^*\{\tilde{v}> t_n\}=\partial^* \{\tilde{u}>t_{n}\}\quad(\operatorname{mod}\;\Hau^1)\qquad \forall n\in\N$$




By  Lemma \ref{regolare},  we can approximate every function $(t_{n+1}-t_n)\mathds{1}_{\{\tilde{v}>t_{n}\}}$ by a
function $\varphi_{n}\in \cont^2(\R^2)$ such that

$$ \mbox{spt}(\varphi_{n})  \subset\!\subset\{\tilde{v}>t_{n}\}$$
$$ \{\tilde{v}>t_{n+1}\}\subset\!\subset \text{spt}(\varphi_{n})$$
\begin{equation}\label{stimafina}
 \|(t_{n+1}-t_n)\mathds{1}_{\{\tilde{v}>t_{n}\}}-\varphi_{n}\|_{\lp{1}}\leq\varepsilon\,2^{-n-2}
\end{equation}
\begin{equation}\label{stimafin}
 \arrowvert G(\Phi_{(t_{n+1}-t_n)\mathds{1}_{\{\tilde{v}>t_{n}\}}})-G(\Phi[\varphi_{n}])\arrowvert\leq
\varepsilon\,2^{-n-2}
\end{equation}
where $\Phi_{(t_{n+1}-t_n)\mathds{1}_{\{\tilde{v}>t_{n}\}}}(t)= \partial \{\tilde{u}>t_{n}\}$ for every $t\in [0, (t_{n+1}-t_n)]$, $\emptyset$ otherwise.\\
\par Finally we  define 
$$v^N(x)=\sum_{n=0}^N \varphi_{n} .$$
that is in $\cont^2(\R^2)$. By \eqref{stimafina} and the Dominated Convergence Theorem, for $N$ large enough
\begin{equation}\label{l4}
 \|\tilde{v}-v^N\|_{\lp{1}}\leq\varepsilon/2.
\end{equation}
By \eqref{stimafin} and the coarea formula, we also have for $N$ large enough
\begin{equation}\label{w4}
 \arrowvert F(v^N)-G(\Phi_{\tilde{v}})\arrowvert = \arrowvert
G(\Phi[v^N])-G(\Phi_{\tilde{v}})\arrowvert\leq \varepsilon/2.
\end{equation}
Finally, by \eqref{l1}, \eqref{l2}, and \eqref{l4},
$$\|u-v^N\|_{\lp{1}} \leq
\|u-\tilde{u}\|_{\lp{1}}+\|\tilde{u}-\tilde{v}\|_{\lp{1}}+\|\tilde{v}-v^N\|_{\lp{1}}\leq \varepsilon/4+ \varepsilon/4 + \varepsilon/2 = \varepsilon,$$
and by \eqref{w1}, \eqref{w2}, and \eqref{w4}
$$\arrowvert G(\Phi)-F(v^N) \arrowvert\leq \arrowvert G(\Phi)-G(\tilde{\Phi})
\arrowvert +  \arrowvert G(\tilde{\Phi})-  G(\Phi_{\tilde{v}})\arrowvert + \arrowvert G(\Phi_{\tilde{v}})-F(v^N)\arrowvert \leq \varepsilon/4+ \varepsilon/4 + \varepsilon/2 = \varepsilon$$

As a straightforward consequence,  we can build a sequence $\{u_h\}\subseteq\cont^2(\R^2)$ such that
$$u_h \overset{\lp{1}}{\longrightarrow}u\quad \mbox{and} \quad G(\Phi)=\underset{h
\rightarrow \infty}{\lim} F(u_h)$$
which implies \eqref{maggiorazione} and the theorem ensues.
\qed

\subsection{Connections between $\oF$ and $\oW$}\label{oFoW}
We further investigate in this section the properties of functions with finite relaxed energy by collecting a few facts about the connection between the energy of a function and the relaxation of the $p$-elastica energy of its level sets.

\par The next proposition generalizes Lemma~\ref{regolare}.

\begin{prop}\label{equiv}
 Let $E\subset \R^2$ be a measurable set such that $\overline{W}(E)<\infty$. Let $u=c\mathds{1}_E$ with $c > 0$. Then 
$$\oF(u)=c \overline{W}(E).$$
\end{prop}

\dimo We first prove that $\oF(u)<\infty$. Since $\overline{W}(E)<\infty$, there exists a sequence of smooth sets $(E_n)$ such that $|E_n\Delta E|\to 0$ and $W(E_n)\to\oW(E)$. Defining $u_n=c\one{E_n}$ it follows from  Lemma~\ref{regolare} that $\oF(u_N)=cW(E_n)$ converges. In addition, $u_n\to u$ in $\lp{1}$ therefore, by the lower semicontinuity of $\oF$, $\oF(u)<\infty$.

Since $\overline{W}(E)<\infty$, by Proposition 6.1 in \cite{BM1}, there exists a finite system of curves $\Gamma$ of class $\wpq{2}{p}$ such that 
\begin{itemize}
\item[(i)]$(\Gamma)\supseteq \partial E$;
\item[(ii)]$|E \Delta \Int(\Gamma)| = 0$;
\end{itemize}
and
$$\overline{W}(E)= \int_{\Gamma} [1+|k_{\Gamma}|^p] \text{d}\mathcal{H}^1.$$
In addition $\Gamma$ minimizes the functional
$$\Gamma \mapsto W(\Gamma) = \int_{\Gamma} [1+|k_{\Gamma}|^p] \text{d}\mathcal{H}^1$$
on the class of all systems of $\wpq{2}{p}$ curves satisfying $(i)$ and $(ii)$. 
Then the function 
$$\Psi : t\in\R \mapsto \Psi(t)=\Gamma\;\text{if $t\in[0,c]$},\;\emptyset\;\text{otherwise}$$
belongs to $\mathscr{A}(u)$ and, for every $\Phi\in \mathscr{A}(u)$,
$$W(\Phi(t))\geq W(\Gamma)\quad \forall t\in [0,c]$$
therefore
$$G(\Phi)\geq G(\Psi) \quad \quad \forall\Phi \in \mathscr{A}(u). $$ 
It follows from Theorem \ref{principale} that
$$\oF(u)= G(\Psi) = c \overline{W}(E).$$
\qed

Using the previous proposition and~\cite{BDP}, we can provide an explicit, and actually trivial, example of a minimizer of $G$ on $\mathscr{A}(u)$.

\begin{ex}{\rm
 Let $u=\mathds{1}_E$ with $E$ like in Figure~\ref{determ}, left, and $L$ the distance between the two cusps. We will prove that  the function
$$\Phi: t\in [0,1]\mapsto \Gamma $$ 
is a minimizer of $G$ on $\mathscr{A}(u)$, being $\Gamma$ the curve in Figure~\ref{determ}, right.

\begin{figure}[h]
\begin{center}
\includegraphics[height=2cm]{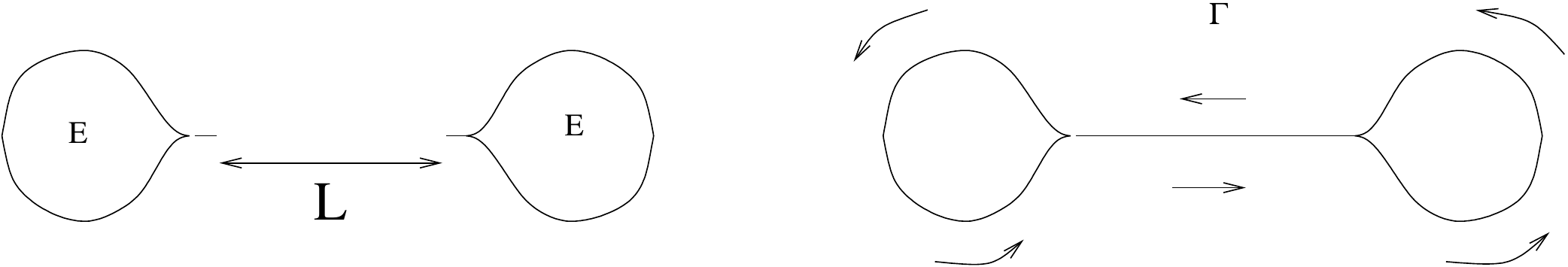}
\caption{The curve $\Gamma$ is canonically associated with the representation of $\oF(\mathds{1}_E)$}\label{determ}
\end{center}
\end{figure}

From the previous proposition we get
\begin{equation}\label{goccia0}
 \oF(u)= \overline{W}(E),
\end{equation}
and by Theorem~8.6 in~\cite{BM1} we have 
\begin{equation}\label{goccia}
\overline{W}(E)= \mathcal{H}^1(\partial E) + \int_{\partial^* E} |k_{\partial^*E}|^p \text{d}\mathcal{H}^1 +2L.
\end{equation}
Now  $\Phi$ belongs to $\mathscr{A}(u)$ and, by definition of $\Gamma$, 
$$G(\Phi)=\mathcal{H}^1(\partial E) + \int_{\partial^* E} |k_{\partial^*E}|^p \text{d}\mathcal{H}^1 +2L$$
therefore
$$\oF(u)=G(\Phi)$$
and, by  Theorem \ref{principale},  $\Phi$ minimizes $G$ on $\mathscr{A}(u)$. }
\end{ex}

The next proposition has been implicitly used in the introduction.
\begin{prop}\label{energiefinie}
 Let $u\in \BV(\R^2)$ with $\oF(u)<\infty$. Then $\oF(u)\geq \ds\int_{\R} \overline{W}(\{u>t\})dt$. In particular, $\overline{W}(\{u>t\})<\infty$ for almost every $t$.
 \end{prop}

\dimo Since $\oF(u)<\infty$ we can find a sequence $\{u_h\}\subseteq \cont^2(\R^2)$ such that
$$\underset{h\rightarrow \infty}{\lim}\,F(u_h) = \oF(u).$$
Then by the coarea formula and Fatou's lemma we get
\begin{equation}\label{inegwill}
 \int_{\R} \underset{h \rightarrow \infty}{\liminf}\; W(\{u_h>t\})\,dt \leq \underset{h \rightarrow \infty}{\lim}\;F(u_h) = \oF(u) <\infty
\end{equation}
and so, for a.e. $t$,
$$\underset{h \rightarrow \infty}{\liminf}\; W(\{u_h>t\}) <\infty.$$
Since for almost every $t$, $|\{u_h>t\} \Delta \{u>t\}|\rightarrow 0$, it follows from the lower semicontinuity of the relaxation that
$$\overline{W}(\{u>t\})\leq \liminf_{h\to\infty}\;W(\{u_h>t\})<\infty\quad{\mbox{for a.e. $t$}}.$$
Therefore, by \eqref{inegwill}, 
$$\oF(u)\geq \int_{\R} \overline{W}(\{u>t\})\mbox{d}t.$$
\qed

\begin{ex}{\rm 
We already mentioned in the introduction and in Remark~\ref{iii-est-nec} the following example of a function $u\in\BV(\R^2)$ such that $\oF(u)> \ds\int_{\R} \overline{W}(\{u>t\})\mbox{d}t$, that we now simply revisit in the perspective of Theorem~\ref{principale}. 

\par Let $u=\mathds{1}_{E\cup F}+\mathds{1}_F$ with $E,F$ as in Figure~\ref{ultima}A) and, in Figure~\ref{ultima}B), the limit systems of curves with their multiplicities corresponding to the independent approximation of $E\cup F$ and $F$, respectively.
\begin{figure}[!htp]
\begin{center}
\includegraphics[height=4cm]{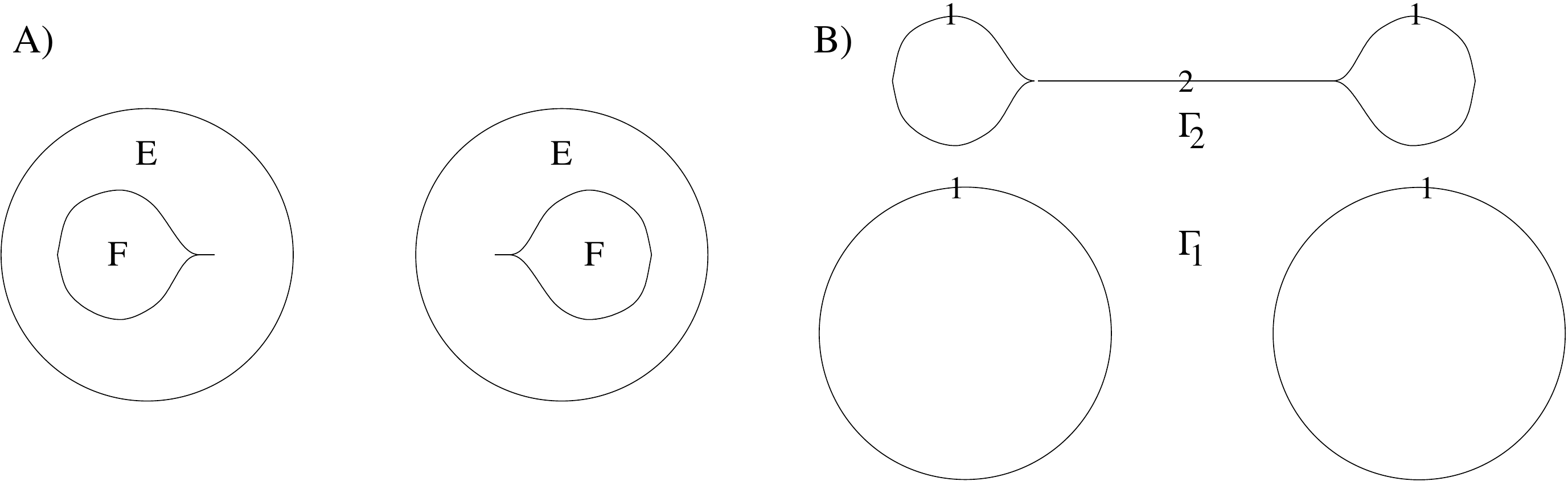}
\caption{If $u=\mathds{1}_{E\cup F}+\mathds{1}_F$, $\oF(u)$ does not coincide with the integral of $\overline{W}(\{u>t\})$.}
\label{ultima}
\end{center}
\end{figure} 

There is no sequence of smooth functions $(u_h)$ approximating $u$ whose level lines concentrate this way because $\Gamma_1$, $\Gamma_2$ together with their densities cannot be contemporaneously approximated (pointwisely) using boundaries of nested sets. In addition, 
$$\int_{\R} \overline{W}(\{u>t\})\mbox{d}t=G(\Phi)$$
where
$$\Phi(t)=\left\{
\begin{array}{ll}
  \Gamma_1  &  \text{if\;} t\in[0,1],\\
  \Gamma_2  &  \text{if\;} t\in]1,2],\\
  \emptyset &  \text{otherwise}.  
\end{array} \right.$$
Obviously, $\Phi\notin \mathscr{A}(u)$, and no system of curves in $\mathscr{A}(u)$ can compete with the energy of $\Phi$ since, to maintain the nesting property, it is necessary to create an additional path between both components of $\Gamma_1$ and the length of this path is at least the distance between the two disks. Therefore,
$$\oF(u) >\int_{\R} \overline{W}(\{u>t\})\mbox{d}t.$$}
\end{ex}

\section{The relaxation problem on a bounded domain}\label{omega}

We consider in this section the generalized elastica functional defined on a bounded domain $\Omega \subset \mathbb{R}^2$ and we compare several definitions for the relaxation problem pointing out their differences. We shall keep the notations $F(u)$ and $\overline{F}(u)$ to denote the generalized elastica energy and its relaxation for a function $u$ defined on $\mathbb{R}^2$.
\par Let $\Omega \subset \mathbb{R}^2$ be an open bounded domain with $\partial \Omega$ Lipschitz. A first definition of a generalized elastica functional on $\Omega$ is:
$$ F_B(\cdot , \Omega): \BV(\Omega) \rightarrow \mathbb{R}$$
$$F_B(u, \Omega)=
\left\{
\begin{array}{ll}
\displaystyle{\int_{\Omega}\left|\nabla u\right|\left(1+\left|\mdiv\frac{\nabla u}{\left|\nabla u\right|}\right|^p\right)\,dx} & \mbox{\,\;if\,} u\in
\cont^2(\Omega)\\
+\infty & \mbox{\,\;otherwise}
\end{array} \right.
$$
with $p>1$, and  $|\nabla u||\mdiv \frac{\nabla u}{\left|\nabla
u\right|}|^p=0$ if $\left|\nabla u\right|=0$.
By definition of the relaxation
$$\overline{F}_B(u, \Omega)= \inf\left\lbrace
\underset{h\rightarrow\infty}{\liminf}\;F_B(u_h, \Omega) : \{u_h\}\subset \cont^2(\Omega), \;
u_h\overset{\lp{1}(\Omega)}{\longrightarrow}u\right\rbrace. $$
Remark first that, by the coarea formula, 
$$F_B(u, \Omega) =  \int_{\mathbb{R}}\int_{\partial\{u>t\}\cap \Omega} (1+ \arrowvert {\kappa}_{\partial\{u>t\}\cap \Omega}\arrowvert^p) \,d\mathcal{H}^1\,dt,\qquad\forall u\in\cont^2(\Omega),$$
so the generalized elastica functional on $\Omega$  depends only on the behavior in $\Omega$ of the level lines of $u$. Therefore, in contrast with what happens in $\R^2$, one cannot restrict to systems of closed curves. Open curves must also be considered, which raises new difficulties, as illustrated in the next example where we exhibit a function with infinite relaxed energy on $\mathbb{R}^2$ and finite relaxed energy on a suitable $\Omega$.

\begin{ex}\label{gocciasola}
{\rm Let $E$, $\Omega$ be the sets drawn in Figure~\ref{gocciao}. Clearly, if $w=\mathds{1}_E$  in $\R^2$, $\oF(v)=\infty$ according to~\cite[Thm 6.4]{BDP}. However, if we consider the relaxation problem on $\Omega$  and  the sequence of  functions $\{w_h\}\subset \cont^2(\Omega)$ having  level lines like the curve $\Gamma$  drawn in Figure~\ref{gocciao}, we get  $\underset{h \rightarrow \infty}{\liminf}\; F_B(w_h, \Omega)<\infty $  thus $\overline{F}_B(w, \Omega)<\infty$.
\begin{figure}[h]
\begin{center}
\includegraphics[height=4.5cm]{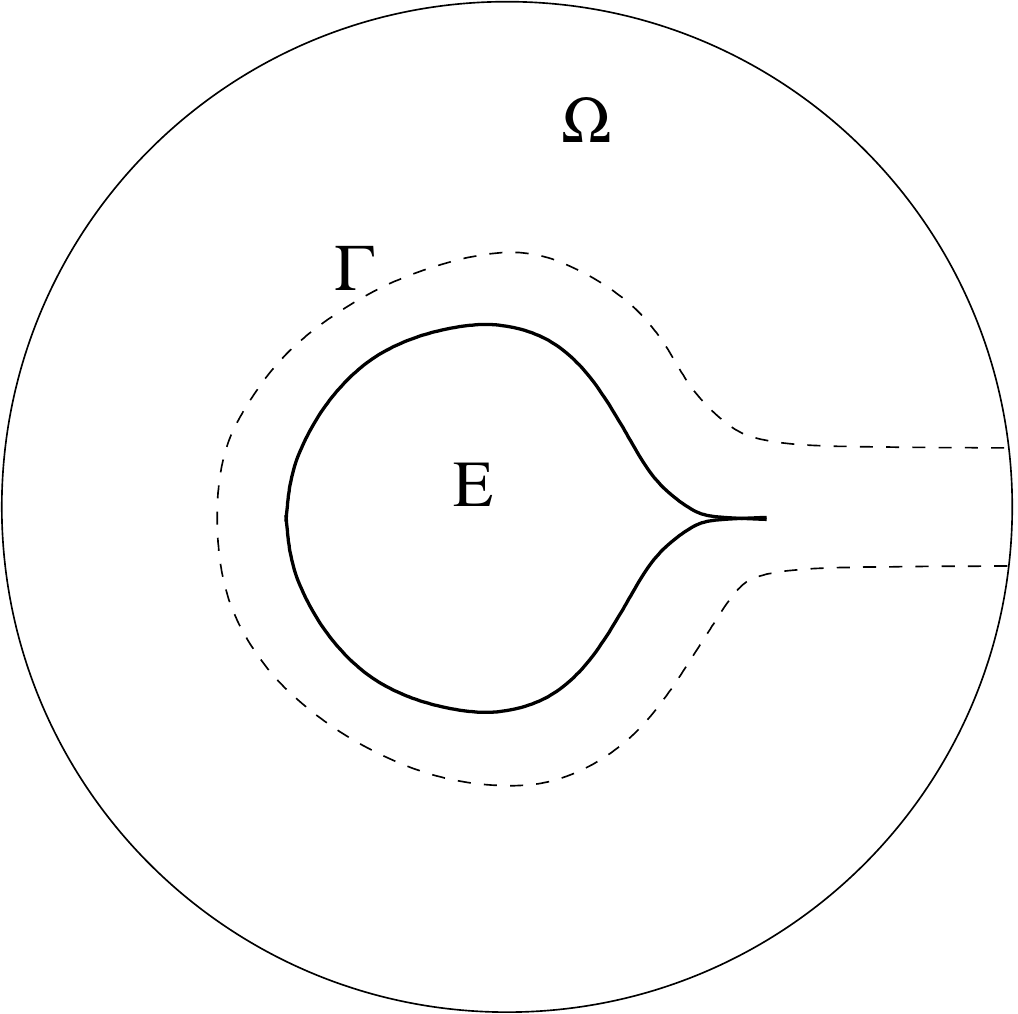}
\caption{$\overline{F}(\mathds{1}_E) = \infty$ but $\overline{F}_B(\mathds{1}_E, \Omega)<\infty$}\label{gocciao}
\end{center}
\end{figure}}
\end{ex}

\par It is a trivial observation that if $u\in\BV(\R^2)$ is such that there exists a sequence $(u_h)$ with $u_h \rightarrow u$ in $\lp{1}(\mathbb{R}^2)$  and $\{F(u_h)\}_h$ is bounded, then, clearly, $\{F_B(u_h|_\Omega, \Omega)\}_h$ is bounded and $\overline{F}_B(u|_\Omega, \Omega)<\infty$. 

Conversely, a natural question is the following: given  $\{u_h\}\subset \cont^2(\Omega)$ a  sequence such that  $F_B(u_h, \Omega)\to\overline{F}_B(u, \Omega)$,
 can we find a sequence $\{v_h\}\subset \cont^2(\mathbb{R}^2)$ with $F(v_h)$ uniformly bounded and $u_h=v_h$ in $\Omega$? In other words, can we say that sequences with bounded energy on $\Omega$ are the restriction to $\Omega$ of sequences with bounded energy on $\mathbb{R}^{2}$?  A positive answer to this question would imply that $\overline{F}_B(\cdot, \Omega)$ coincides in $\BV(\Omega)$ with $L(\cdot, \Omega)$
where 
$$L(u,\Omega) := \inf\left\lbrace
\underset{h\rightarrow\infty}{\liminf}\;F_B(u_h|_\Omega, \Omega) : \{u_h\}\subset \cont^2(\mathbb{R}^2), \;
u_h|_{\Omega}\overset{\lp{1}(\Omega)}{\longrightarrow}u\;,\;\; \underset{h}{\sup}
\;F(u_h)<\infty\right\rbrace, $$
with the convention $\inf\;\emptyset=\infty$.  

\par A simple example of a function with finite $L(\cdot,\Omega)$ energy is the function $w=\one{E}$ of Example~\ref{gocciasola}. Take indeed the image $F$ of $E$ obtained by symmetry with respect to a vertical axis arbitrarily chosen at the right of $\Omega$. Then, $\partial(E\cup F)$ being smooth except at an even number of cusps,  $\oF(E\cup F)<+\infty$ according to Theorem~6.3 in~\cite{BDP}. Thus there exists a sequence of smooth functions $(u_h)$ that approximates $\one{E\cup F}$ in $\R^2$, belongs to the set $\left\lbrace
\{u_h\}\subset \cont^2(\mathbb{R}^2), \;
u_h|_{\Omega}\overset{\lp{1}(\Omega)}{\longrightarrow}u\;,\;\; \underset{h}{\sup}
\;F(u_h)<\infty\right\rbrace$ and is such that $\sup_h F_B(u_h|_\Omega, \Omega)<\infty$, therefore $L(u,\Omega)<\infty$.

\par The answer to the question above is negative in general as shown by the next example.

\begin{ex}\label{gocciadoppia}
{\rm Let $E, \Omega$ be the sets drawn in Figure~\ref{gocciad} and  let $u=\mathds{1}_E$. We know from Theorem 4.1 in~\cite{BDP} that any set $E$ with finite relaxed energy is such that $\partial* E$ has a continuous unoriented tangent, which is obviously not the case here, thus $L(\mathds{1}_E, \Omega) = \infty$. Roughly speaking, every sequence  $\{v_h\}\subset \cont^2(\mathbb{R}^2)$ converging to $u$ in $\Omega$ has to approximate the angle in $p\in\partial\Omega$ formed by the two cusps and so $F(v_h)\rightarrow \infty $. Now, if we consider the relaxation problem on $\Omega$  and  the sequence of  functions $\{u_h\}\subset \cont^2(\Omega)$ having  level lines like the curve $\Gamma$  drawn in Figure~\ref{gocciad}, there is no singularity in $\Omega$ and $\underset{h \rightarrow \infty}{\liminf}\; F_B(u_h, \Omega)<\infty $  thus $\overline{F}_B(u, \Omega)<\infty$.
\begin{figure}[h]
\begin{center}
\includegraphics[height=4.5cm]{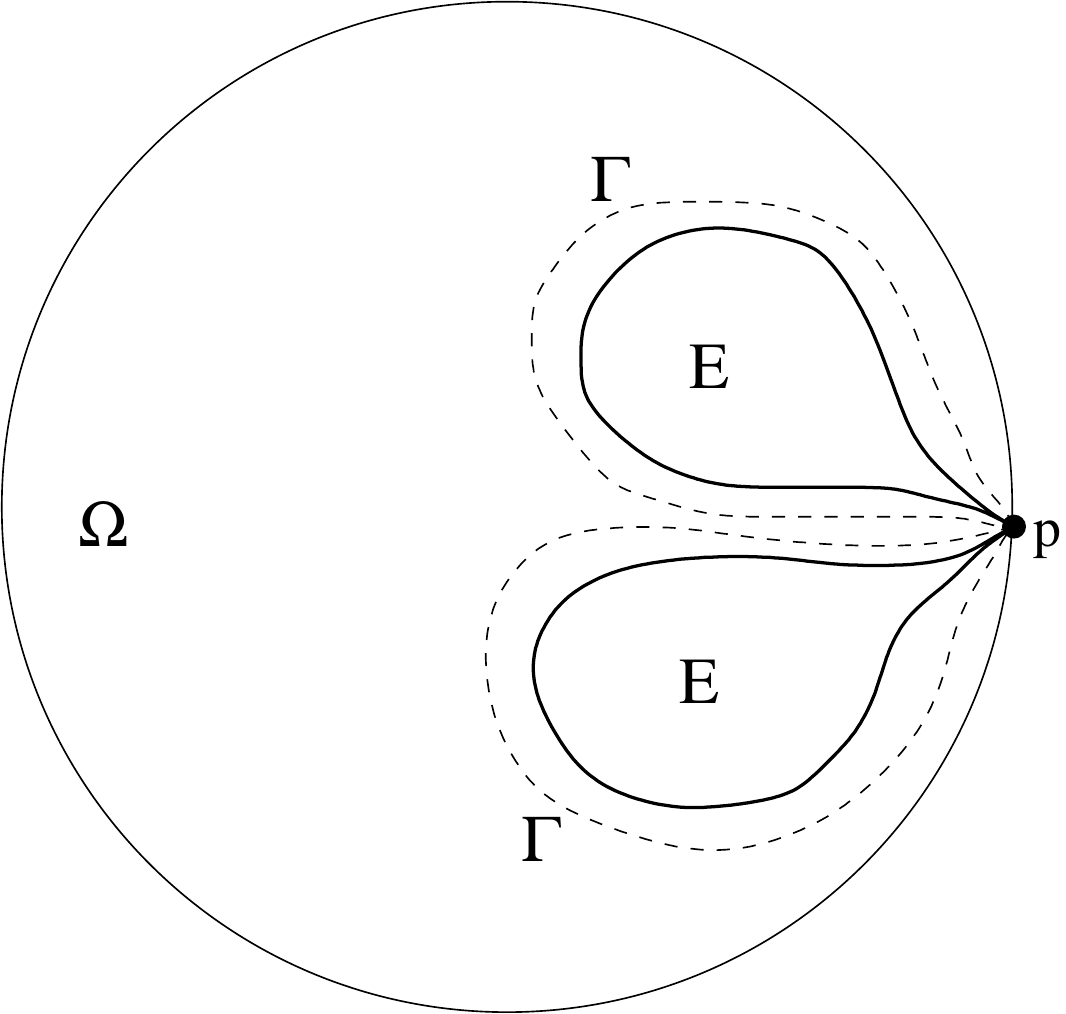}
\caption{$L(\mathds{1}_E, \Omega) = \infty$ but $\overline{F}_B(\mathds{1}_E, \Omega)<\infty$}\label{gocciad}
\end{center}
\end{figure}}
\end{ex}

Another possible way to define the generalized elastica functional on an open bounded domain with Lipschitz boundary is the following, where $\BV_0(\Omega)$ denotes the space of functions of bounded variation defined in $\Omega$ with null trace on $\partial\Omega$:
$$F_B^0(\cdot, \Omega): \BV_0(\Omega) \rightarrow \mathbb{R}$$
$$F_B^0(u,\Omega)=
\left\{
\begin{array}{ll}
\displaystyle{\int_{\Omega}\left|\nabla u\right|\left(1+\left|\mdiv
\frac{\nabla u}{\left|\nabla u\right|}\right|^p\right)\,dx}  &\mbox{\,\;if\,} u\in
\cont^2_c(\Omega)\\
+\infty & \mbox{\,\;otherwise}
\end{array} \right.$$
with $p>1$, and  $\left|\nabla u\right|\left(1+\left|\mdiv \frac{\nabla u}{\left|\nabla
u\right|}\right|^p\right)=0$ if $\left|\nabla u\right|=0$.

\par In this case the relaxation in $\BV_0(\Omega)$ is defined by
$$\overline{F}_B^0(u,\Omega)=\inf\left\lbrace
\underset{h\rightarrow\infty}{\liminf}\;F_B^0(u_h, \Omega) : \{u_h\}\subset \cont^2_c(\Omega), \;
u_h\overset{\lp{1}(\Omega)}{\longrightarrow}u\right\rbrace.$$
Remark that the function of Example~\ref{gocciasola} is in $\BV_0(\Omega)$ and has infinite $F_B^0$ energy. In contrast, the function of Example~\ref{gocciadoppia} is not in the domain of $\overline{F}_B^0$ since it is not in $\BV_0(\Omega)$. It would not make sense to extend, by approximation, the definition of $\overline{F}_B^0$ to the functions of $\BV(\Omega)\setminus\BV_0(\Omega)$ since $\one{\Omega}$ is such function, has no level line in $\Omega$, but 
$$\inf\left\lbrace
\underset{h\rightarrow\infty}{\liminf}\;F_B^0(u_h, \Omega) : \{u_h\}\subset \cont^2_c(\Omega), \;
u_h\overset{\lp{1}(\Omega)}{\longrightarrow}\one{\Omega}\right\rbrace=W(\Omega)>0.$$

\par The next proposition states the very natural localization property that a function with compact support that has finite energy can be approximated, in $\lp{1}$ and in energy, by a sequence of functions with compact support. The approximating sequence is built directly from the collection of curves that cover the level lines of the function.
\begin{prop}\label{legameY} 
If $u\in\BV(\mathbb{R}^2)$ has compact support and $\overline{F}(u) < \infty$, there exists an open and bounded domain $\Omega$ that contains $\text{spt}\, u$ and satisfies 
 $$\overline{F}(u)=\overline{F}_B^0(u|_\Omega,\Omega).$$
\end{prop}
\dimo Using Theorem \ref{principale} we have $\overline{F}(u)= G(\Phi)$ where $\Phi$ is a minimizer of $G$ on  $\mathscr{A}(u)$. Since $u$ has compact support on $\mathbb{R}^2$ , the definition of $\mathscr{A}(u)$ and $G(\Phi)<\infty$ imply the existence of a bounded and open domain $\Omega$ such that
\begin{equation}\label{inclusion}
\bigcup_{t\in \mathbb{R}} (\Phi(t)) \subset \subset \Omega.
\end{equation}
\par Since $\overline{F}(u) \leq \inf\left\lbrace
\underset{h\rightarrow\infty}{\liminf}\;F_B^0(u_h, \Omega) : \{u_h\}\in \cont^2_c(\Omega), \; 
u_h\overset{\lp{1}(\Omega)}{\longrightarrow}u\right\rbrace$, using~\eqref{inclusion} and the same approximation arguments developed in the proof of Theorem~\ref{principale} together with the fact that $\{u>0\}$ is relatively compactly, we can define a sequence $\{u_h\}\in \cont^2_c(\Omega)$ such that
$u_h\overset{\lp{1}(\Omega)}{\longrightarrow}u$ and $F(u_h)=F_B^0(u_h, \Omega) \rightarrow G(\Phi)$. Then 
$$G(\Phi)= \overline{F}(u) \geq \inf\left\lbrace
\underset{h\rightarrow\infty}{\liminf}\;F_B^0(u_h, \Omega) : \{u_h\}\in \cont^2_c(\Omega), \; 
u_h\overset{\lp{1}(\Omega)}{\longrightarrow}u\right\rbrace$$
and the proposition ensues.
\qed


\par The next example illustrates the difference between $\overline{F}_B(u,\Omega)$ and $\overline{F}_B^0(u,\Omega)$, that are defined using different function spaces for the approximation ($\cont^2(\Omega)$ for the former and $\cont^2_c(\Omega)$ for the latter).

\begin{ex}
{\rm Let $E, \Omega$ be the sets drawn in Figure~\ref{ultimaborn} and let  $u=\mathds{1}_E$. 
\begin{figure}[h]
\begin{center}
\includegraphics[height=4cm]{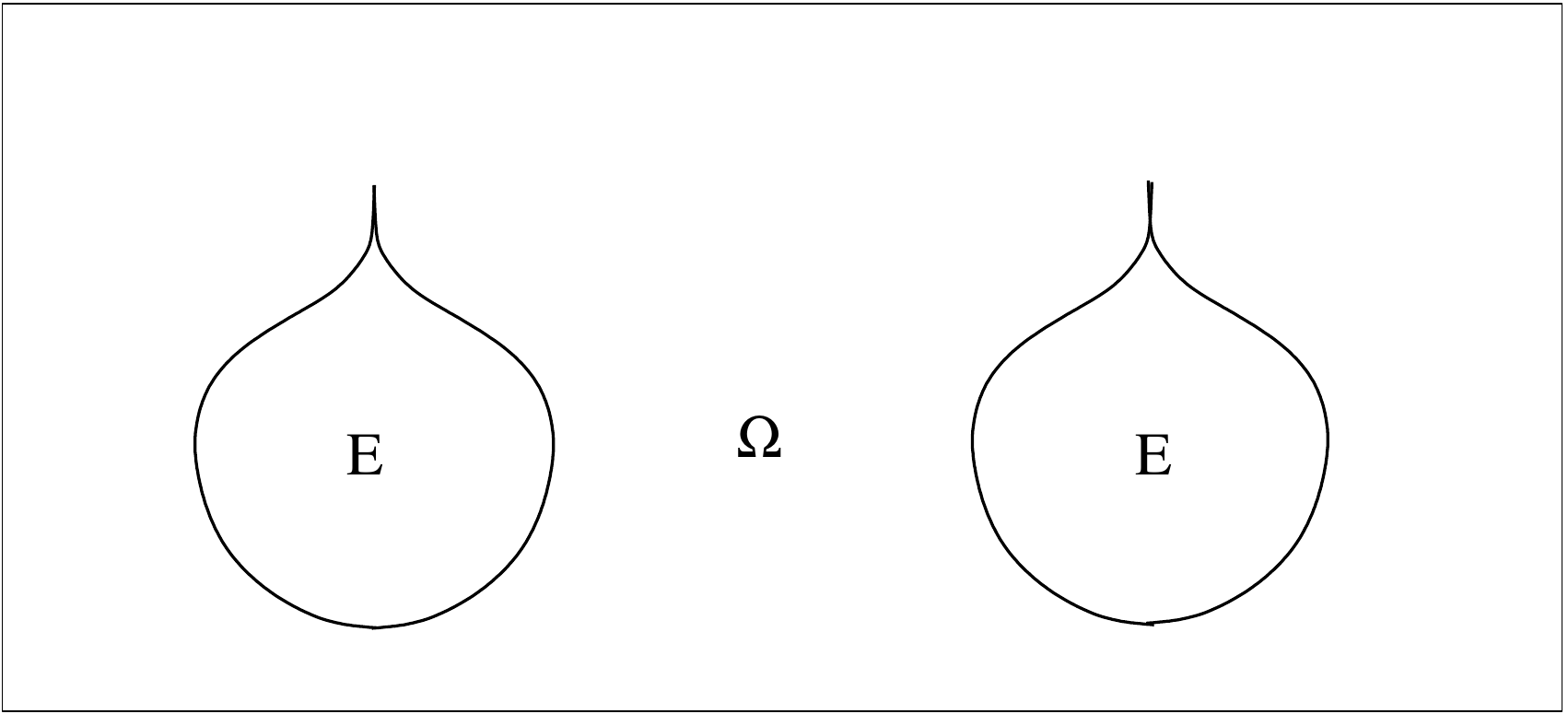}
\caption{$u=\mathds{1}_E\in BV(\Omega)$}\label{ultimaborn}
\end{center}
\end{figure}

From the properties of the relaxation, from the previous proposition, and thanks to Theorem~8.6 in~\cite{BM1}, we have 
$$\overline{F}(u)=\inf\left\lbrace
\underset{h\rightarrow\infty}{\liminf}\;F_B^0(u_h, \Omega) : \{u_h\}\in \cont^2_c(\Omega), \; 
u_h\overset{\lp{1}(\Omega)}{\longrightarrow}u\right\rbrace = G(\Phi)$$
where $\Omega$  is given by the previous proposition and $\Phi(t)$ looks for every $t\in[0;1]$ like the curve  $\gamma_1$ drawn with its multiplicity in Figure~\ref{gamma12}A.

\begin{figure}[h]
\begin{center}
\includegraphics[height=4cm]{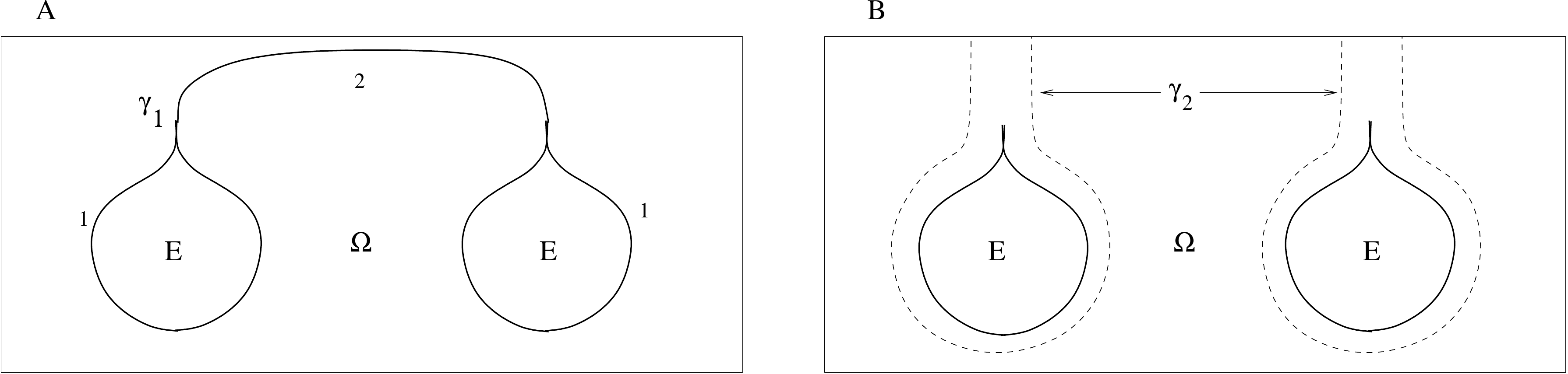}
\caption{Relaxation using sequences belonging to $\cont^2_c(\Omega)$ (left) or $\cont^2(\Omega)$ (right)}\label{gamma12}
\end{center}
\end{figure}

\par However, if we consider the sequence  of  functions $\{u_h\}\subset \cont^2(\Omega)$ having  level lines like the line $\gamma_2$  drawn in Figure~\ref{gamma12}B we get  $\underset{h \rightarrow \infty}{\liminf}\; F_B(u_h, \Omega)< G(\Phi) $. Therefore
$$ \overline{F}_B(u, \Omega) < \overline{F}_B^0(u_h, \Omega).$$ }
\end{ex}
\subsection*{Acknowledgements}
We warmly thank Vicent Caselles and Matteo Novaga for several discussions during the elaboration of this paper.\\
This work was supported by the French "Agence Nationale de la
Recherche" (ANR), under grant FREEDOM (ANR07-JCJC-0048-01).
\bibliographystyle{plain}
\bibliography{biblio2}

\begin{thebibliography}{10}

\bibitem{AFP}
L.~Ambrosio, N.~Fusco, and D.~Pallara.
\newblock {\em Functions of Bounded Variation and Free Discontinuity Problems}.
\newblock Oxford Science Publications, 2000.

\bibitem{AGP}
L.~Ambrosio, M.~Gobbino, and D.~Pallara.
\newblock Approximation problems for curvature varifolds.
\newblock {\em J. Geometric Analysis}, 8(1):1--19, 1998.

\bibitem{AM}
L.~Ambrosio and S.~Masnou.
\newblock A direct variational approach to a problem arising in image
  reconstruction.
\newblock {\em Interfaces and Free Boundaries}, 5:63--81, 2003.

\bibitem{Ballester}
C.~Ballester, M.~Bertalmio, V.~Caselles, G.~Sapiro, and J.~Verdera.
\newblock Filling-in by joint interpolation of vector fields and gray levels.
\newblock {\em IEEE Trans. On Image Processing}, 10(8):1200--1211, 2001.

\bibitem{BDP}
G.~Bellettini, G.~Dal Maso, and M.~Paolini.
\newblock Semicontinuity and relaxation properties of a curvature depending
  functional in 2{D}.
\newblock {\em Annali della Scuola Normale di Pisa, Classe di Scienze, $4^e$
  s\'erie}, 20(2):247--297, 1993.

\bibitem{BM1}
G.~Bellettini and L.~Mugnai.
\newblock Characterization and representation of the lower semicontinuous
  envelope of the elastica functional.
\newblock {\em Ann. Inst. H. Poincar\'e, Anal. non Lin\'eaire}, 21(6):839--880,
  2004.

\bibitem{BM}
G.~Bellettini and L.~Mugnai.
\newblock A varifold representation of the relaxed elastica functional.
\newblock {\em Journal of Convex Analysis}, 14(3):543--564, 2007.

\bibitem{Br}
K.A. Brakke.
\newblock The motion of a surface by its mean curvature.
\newblock {\em Mathematical notes (20), Princeton University Press}, 1978.

\bibitem{CGMP}
F.~Cao, Y.~Gousseau, S.~Masnou, and P.~P\'{e}rez.
\newblock Geometrically guided exemplar-based inpainting.
\newblock {\em SIAM Journal on Imaging Sciences}, 2011.

\bibitem{ChanKangShen}
T.F. Chan, S.H. Kang, and J.~Shen.
\newblock Euler's elastica and curvature based inpainting.
\newblock {\em SIAM Journal of Applied Math.}, 63(2):564--592, 2002.

\bibitem{CittiSarti}
G.~Citti and A.~Sarti.
\newblock A cortical based model of perceptual completion in the
  roto-translation space.
\newblock {\em J. Math. Imaging Vis.}, 24(3):307--326, 2006.

\bibitem{DalMaso}
G.~Dal~Maso.
\newblock {\em An introduction to $\Gamma$-convergence}, volume~8 of {\em
  Progress in Nonlinear Diff. Equ. and their Appl.}
\newblock Birkha{\"u}ser, Boston, 1993.

\bibitem{GT}
D.~Gilbarg and N.S. Trudinger.
\newblock {\em Elliptic Partial Differential Equations of Second Order}.
\newblock Springer, 1998.

\bibitem{G}
A.~Gray.
\newblock {\em Modern differential geometry of curves and surfaces with
  Mathematica}.
\newblock CRC Press, 2000.

\bibitem{LM}
G.P. Leonardi and S.~Masnou.
\newblock Locality of the mean curvature of rectifiable varifolds.
\newblock {\em Adv. Calc. of Var.}, 2(1):17--42, 2009.

\bibitem{MasnouMorel}
S.~Masnou and J.-M. Morel.
\newblock Level lines based disocclusion.
\newblock In {\em 5th IEEE Int. Conf. on Image Processing, Chicago, Illinois,
  October 4-7}, 1998.

\bibitem{MM}
S.~Masnou and J.M. Morel.
\newblock On a variational theory of image amodal completion.
\newblock {\em Rendiconti del Seminario Matematico della Universit\`{a} di
  Padova}, 116:211--252, 2006.

\bibitem{MasnouNardi2}
S.~Masnou and G.~Nardi.
\newblock Gradient {Y}oung measures, varifolds, and a generalized {W}illmore
  functional.
\newblock submitted, 2011.

\bibitem{Menne}
U.~Menne.
\newblock Second order rectifiability of integral varifolds of locally bounded
  first variation.
\newblock arXiv:0808.3665v3 [math.DG]., 2010.

\bibitem{Petitot}
J.~Petitot.
\newblock Neurogeometry of {V}1 and {K}anizsa contours.
\newblock {\em Axiomathes}, 13:347--363, 2003.

\bibitem{RY}
D.~Revuz and M.~Yor.
\newblock {\em Continuous martingales and Brownian Motion}.
\newblock Springer, 1998.

\bibitem{rudin-real-complex}
W.~Rudin.
\newblock {\em Real and complex analysis}.
\newblock McGraw-Hill Book Co., New York, 3rd edition, 1987.

\bibitem{savare}
G.~Savar{\'e}.
\newblock Personal communication.

\bibitem{S}
R.~Sch{\"a}tzle.
\newblock Lower semicontinuity of the {W}illmore functional for currents.
\newblock {\em J. Differential Geometry}, 8:437--456, 2009.

\bibitem{Si}
L.~Simon.
\newblock Lectures on geometric measure theory.
\newblock volume~3 of {\em Proc. Centre for Math. Analysis,}. Australian Nat.
  Univ., 1983.

\end{thebibliography}
\end{document}